\documentclass[12pt,reqno]{amsart}
\usepackage{fullpage,mathpple,mathrsfs}
\newtheorem{theorem}{Theorem}
\newtheorem{lemma}[theorem]{Lemma}
\newtheorem{cor}[theorem]{Corollary}
\newtheorem{conjecture}[theorem]{Conjecture}
\newtheorem{prop}[theorem]{Proposition}
\newtheorem{definition}{Definition}[theorem]

\renewcommand{\mod}[1]{{\ifmmode\text{\rm\ (mod~$#1$)}\else\discretionary{}{}{\hbox{ }}\rm(mod~$#1$)\fi}}
\newcommand{\Sk}[2]{S_k(\Gamma_#1(#2))}
\newcommand{\Skstar}[2]{S_k^*(\Gamma_#1(#2))}
\newcommand{\Skplus}[2]{S_k^+(\Gamma_#1(#2))}

\newcommand{\Z}{{\mathbb Z}}
\newcommand{\N}{{\mathbb N}}
\newcommand{\C}{{\mathbb C}}
\newcommand{\ep}{\varepsilon}
\newcommand{\exdiv}{\mathrel{\|}}

\begin{document}

\title{Dimensions of the Spaces of Cusp Forms and Newforms on $\Gamma_0(N)$ and $\Gamma_1(N)$}
\author{Greg Martin}
\address{Department of Mathematics\\University of British Columbia\\Room 121, 1984 Mathematics Road\\Vancouver, BC V6T 1Z2}
\email{gerg@math.ubc.ca}
\maketitle

\section{Introduction}

The study of modular forms on congruence groups was initiated by Hecke and Petersson in the 1930s and, at least when the weight $k$ is an integer exceeding~1, is quite well understood. In particular, formulas for the dimensions of the spaces of modular forms and cusp forms on the congruence groups
\[
\Gamma_0(N) = \Big\{ \Big( \!\begin{array}{c@{\ }c}a&b\\c&d\end{array} \!\Big)\! \in SL_2(\Z) \colon c\equiv0\mod N \Big\}
\]
and
\[
\Gamma_1(N) = \Big\{ \Big( \!\begin{array}{c@{\ }c}a&b\\c&d\end{array} \!\Big)\! \in SL_2(\Z) \colon a\equiv d\equiv1\mod N,\, c\equiv0\mod N \Big\}
\]
are known \cite{miyake,Shimura} (see Propositions \ref{g0.formula.prop} and \ref{g1.formula.prop} below). The structure of these spaces of cusp forms was clarified by the work of Atkin and Lehner \cite{AL}, who exhibited their orthogonal decomposition with respect to the Petersson inner product into spaces of cuspidal newforms. Until now, however, the dimensions of the spaces of newforms
could only be calculated recursively (in terms of the corresponding dimensions for divisors of the level $N$) and thus were rather poorly understood in general.

In this paper we present closed formulas for the dimensions of the spaces of weight-$k$ cuspidal newforms on $\Gamma_0(N)$ and $\Gamma_1(N)$, for all integers $k\ge2$. The formulas consist of linear combinations of multiplicative functions of $N$, with coefficients depending on $k$; in particular, they have the same level of simplicity as the formulas for the dimensions of the full spaces of cusp forms on these modular groups. As an application of the new formulas, we derive simple upper and lower bounds for the dimensions of these spaces of newforms for all $k\ge2$. We also calculate all positive integers $N$ for which the dimension of the space of newforms of weight 2 on $\Gamma_0(N)$ is at most 100, and we prove the validity of certain inequalities and identities for these dimensions observed empirically by Bennett. Finally, the question of the dimensions of these spaces on average over $N$ does not seem to have been raised even for the full spaces of cusp forms. We calculate the average orders both of the dimensions of the spaces of weight-$k$ cusp forms on $\Gamma_0(N)$ and $\Gamma_1(N)$ and of the dimensions of the subspaces of newforms. In addition, we establish analogues of all these results for the numbers of nonisomorphic automorphic representations associated with these spaces of modular forms.

We now set some notation with which to describe our results. Let $\Sk0N$ denote the space of cusp forms on $\Gamma_0(N)$ of weight $k$ and $\Skplus0N$ the space of newforms on $\Gamma_0(N)$ of weight $k$. Let $g_0(N,k)$ and $g_0^+(N,k)$ denote the dimensions of $\Sk0N$ and  $\Skplus0N$, respectively. Our formula for $g_0^+(N,k)$ involves several multiplicative functions that we shall define shortly. Recall that a function $f$, not identically zero, is multiplicative if $f(mn)=f(m)f(n)$ whenever $m$ and $n$ are relatively prime. It follows that $f(1)=1$ and that $f$ is completely determined by its values on prime powers. Some common examples of multiplicative functions that will be useful to us are Euler's totient function $\phi(n)$ and the M\"obius function $\mu(n)$; also $\omega(n)$, the number of distinct prime factors of $n$, and $\tau(n)$, the number of positive divisors of $n$; and finally the delta function at 1,
\begin{equation}
\delta(n) = \begin{cases}
1, &\text{if }n=1, \\ 0, &\text{otherwise.}
\end{cases}
\label{delta.def}
\end{equation}
Our first theorem shows that $g_0^+(N,k)$ can be expressed as a linear combination of multiplicative functions of $N$, with the coefficients depending on $k$.

%:Theorem - formula for g_0^+
\begin{theorem}
For any even integer $k\ge2$ and any integer $N\ge1$, we have
\[
g_0^+(N,k) = \tfrac{k-1}{12}Ns_0^+(N) - \tfrac12\nu^+_\infty(N)+c_2(k) \nu^+_2(N) +c_3(k) \nu^+_3(N) + \delta\big( \tfrac k2 \big)\mu(N),
\]
where the functions $s_0^+$, $\nu_\infty^+$, $\nu_2^+$, $\nu_3^+$, $c_2$, and $c_3$ are defined in Definitions \ref{s0plusdef}--\ref{c3def} below.
\label{g0plus.formula.theorem}
\end{theorem}

%:(Definitions needed to state the g_0^+ formula)
We remark that the restriction that $k$ be even is natural, since there are no modular forms on $\Gamma_0(N)$ of odd integer weight, that is, $g_0(N,k)=0$ and hence $g_0^+(N,k)=0$ when $k$ is odd. We promptly give the definitions of the six functions in the statement of Theorem~\ref{g0plus.formula.theorem}. In the definitions of the multiplicative functions and throughout this paper, $p$ always denotes a prime number.

\begin{definition}
$s_0^+$ is the multiplicative function satisfying

\centerline{$
s_0^+(p) = 1-\tfrac1p,\;
s_0^+(p^2) = 1-\tfrac1p-\tfrac1{p^2},\;
\text{and } s_0^+(p^\alpha) = \big(1-\tfrac1p\big)\big(1-\tfrac1{p^2}\big) \text{ for } \alpha\ge3.
$}
\label{s0plusdef}
\end{definition}

\begin{definition}
$\nu_\infty^+$ is the multiplicative function satisfying

\centerline{$
\nu^+_\infty(p^\alpha)=0 \text{ for $\alpha$ odd},\;
\nu^+_\infty(p^2)=p-2,\;
\text{and } \nu^+_\infty(p^\alpha) = p^{\alpha/2-2}(p-1)^2 \text{ for $\alpha\ge4$ even}.
$}
\label{nuinftyplusdef}
\end{definition}

\begin{definition}
$\nu_2^+$ is the multiplicative function satisfying:
\begin{itemize}
\item $\nu^+_2(2)=-1$,\, $\nu^+_2(4)=-1$,\, $\nu^+_2(8)=1$,\, and $\nu^+_2(2^\alpha)=0$ for $\alpha\ge4$;
\item if $p\equiv1\mod4$ then $\nu^+_2(p)=0$,\, $\nu^+_2(p^2)=-1$,\, and $\nu^+_2(p^\alpha)=0$ for $\alpha\ge3$;
\item if $p\equiv3\mod4$ then $\nu^+_2(p)=-2$,\, $\nu^+_2(p^2)=1$,\, and $\nu^+_2(p^\alpha)=0$ for $\alpha\ge3$.
\end{itemize}
\label{nu2plusdef}
\end{definition}

\begin{definition}
$\nu_3^+$ is the multiplicative function satisfying:
\begin{itemize}
\item $\nu^+_3(3)=-1$,\, $\nu^+_3(9)=-1$,\, $\nu^+_3(27)=1$,\, and $\nu^+_3(3^\alpha)=0$ for $\alpha\ge4$;
\item if $p\equiv1\mod3$ then $\nu^+_3(p)=0$,\, $\nu^+_3(p^2)=-1$,\, and $\nu^+_3(p^\alpha)=0$ for $\alpha\ge3$;
\item if $p\equiv2\mod3$ then $\nu^+_3(p)=-2$,\, $\nu^+_3(p^2)=1$,\, and $\nu^+_3(p^\alpha)=0$ for $\alpha\ge3$.
\end{itemize}
\label{nu3plusdef}
\end{definition}

\begin{definition}
$c_2$ is the function defined by $c_2(k) = \frac14 + \big\lfloor\frac k4\big\rfloor - \frac k4$.
\label{c2def}
\end{definition}

\begin{definition}
$c_3$ is the function defined by $c_3(k) = \frac13 + \big\lfloor\frac k3\big\rfloor - \frac k3$.
\label{c3def}
\end{definition}

We remark that as this manuscript was being prepared, a paper of Halberstadt and Kraus \cite{HK} appeared, in the appendix of which they independently established the special case of Theorem \ref{g0plus.formula.theorem} where $k=2$.

The formula given in Theorem \ref{g0plus.formula.theorem} provides a method of computing $g_0^+(N,k)$ that is much faster than the recursive formula \eqref{g0plus.in.terms.of.g0} below.
% For example, as mentioned earlier, it is a matter of mere minutes to use this formula to compute $g_0^+(N,2)$ for all $N$ up to 100,000.
In Section \ref{up.to.100.section} we show how to use such a computation to determine the complete list of positive integers $N$ such that $g_0^+(N,2)$ is at most 100. Previously, exhaustive lists of those $N$ for which $g_0^+(N,2)=j$ had been given \cite{HK} only for $j=0,1,2,3$. We also gather evidence supporting the assertion that every nonnegative integer is a value of the function $g_0^+(N,2)$, but we refute this assertion for $g_0(N,2)$ itself---the first omitted value is 150.

Moreover, the formula in Theorem \ref{g0plus.formula.theorem} is much more amenable to analysis of the behavior of the function $g_0^+(N,k)$. For example, the coefficients of the last four multiplicative functions $\nu_\infty^+$, $\nu_2^+$, $\nu_3^+$, and $\mu$ in this formula are all bounded functions of $k$. Therefore we can immediately conclude that when $N$ is fixed, the dimension $g_0^+(N,k)$ grows roughly linearly with $k$; more precisely,
\[
g_0^+(N,k) = \tfrac{Ns_0^+(N)}{12}k + O_N(1).
\]
Two further concrete examples of the usefulness of the explicit formula in Theorem \ref{g0plus.formula.theorem} are provided by the following two results. These theorems establish the validity of certain identities and inequalities proposed by Bennett (personal communication) on the basis of numerical observations.

%:Theorem - Bennett's bound
\begin{theorem}
For all positive integers $N$, we have $g_0^+(N,2) \le (N+1)/12$, with equality holding if and only if either $N=35$ or $N$ is a prime that is congruent to $11\mod{12}$.
\label{bennett.bound.theorem}
\end{theorem}

%:Theorem - Bennett's inequality
\begin{theorem}
Let $N\ge3$ be an odd squarefree integer. Then $g_0^+(2^\alpha N,k) = (k-1)2^{\alpha-5}\phi(N)$ for every integer $\alpha\ge4$; in particular, $g_0^+(32N,k) = (k-1)\phi(N)$. In addition, we have $g_0^+(2N,k) \le (k-1)\phi(N)$.
\label{bennett.inequality.theorem}
\end{theorem}

The method of proof of Theorem \ref{g0plus.formula.theorem} can also be used to establish a similar formula for the number of nonisomorphic automorphic representations associated with $\Sk0N$, which we denote by $g_0^*(N,k)$. (See the proof of Theorem \ref{g0star.formula.theorem} in Section \ref{main.theorems.proofs.section} for a more precise definition of the number in question.) Our next theorem shows that $g_0^*(N,k)$ can also be expressed as a linear combination of multiplicative functions of $N$.

%:Theorem - formula for g_0^*
\begin{theorem}
For any even integer $k\ge2$ and any integer $N\ge1$, we have
\[
g_0^*(N,k) = \tfrac{k-1}{12}Ns_0^*(N) - \tfrac12\nu^*_\infty(N)+c_2(k) \nu^*_2(N) +c_3(k) \nu^*_3(N) + \delta\big( \tfrac k2 \big)\delta(N),
\]
where the functions $c_2$, $c_3$, $s_0^*$, $\nu_\infty^*$, $\nu_2^*$, and $\nu_3^*$ are defined in Definitions \ref{c2def}--\ref{c3def} above and Definitions \ref{s0stardef}--\ref{nu3stardef} below.
\label{g0star.formula.theorem}
\end{theorem}

%:(Definitions needed to state the g_0^* formula)
The definitions of the four new functions in the statement of Theorem \ref{g0star.formula.theorem} are as follows.

\begin{definition}
$s_0^*$ is the multiplicative function satisfying

\centerline{$
s_0^*(p) = 1\;
\text{and } s_0^*(p^\alpha) = 1-\tfrac1{p^2} \text{ for } \alpha\ge2.
$}
\label{s0stardef}
\end{definition}

\begin{definition}
$\nu_\infty^*$ is the multiplicative function satisfying

\centerline{$
\nu^*_\infty(p)=1\;
\text{and } \nu^*_\infty(p^\alpha) = p^{\lfloor\alpha/2-1\rfloor}(p-1) \text{ for $\alpha\ge2$}.
$}
\label{nuinftystardef}
\end{definition}

\begin{definition}
$\nu_2^*$ is the multiplicative function satisfying:
\begin{itemize}
\item $\nu^*_2(2)=0$,\, $\nu^*_2(4)=-1$,\, and $\nu^*_2(2^\alpha)=0$ for $\alpha\ge3$;
\item if $p\equiv1\mod4$ then $\nu^*_2(p)=1$\, and $\nu^*_2(p^\alpha)=0$ for $\alpha\ge2$;
\item if $p\equiv3\mod4$ then $\nu^*_2(p)=-1$\, and $\nu^*_2(p^\alpha)=0$ for $\alpha\ge2$.
\end{itemize}
\label{nu2stardef}
\end{definition}

\begin{definition}
$\nu_3^*$ is the multiplicative function satisfying:
\begin{itemize}
\item $\nu^*_3(3)=0$,\, $\nu^*_3(9)=-1$,\, and $\nu^*_3(3^\alpha)=0$ for $\alpha\ge3$;
\item if $p\equiv1\mod3$ then $\nu^*_3(p)=1$\, and $\nu^*_3(p^\alpha)=0$ for $\alpha\ge2$;
\item if $p\equiv2\mod3$ then $\nu^*_3(p)=-1$\, and $\nu^*_3(p^\alpha)=0$ for $\alpha\ge2$.
\end{itemize}
\label{nu3stardef}
\end{definition}

Theorem \ref{g0star.formula.theorem} allows a very short proof of a result of Gekeler \cite{gekeler} in the case where the level $N$ is squarefree:

%:Corollary - when N is squarefree
\begin{cor}
Let $k\ge2$ be an even integer, and let $N\ge1$ be a squarefree integer, with $N>1$ if $k=2$. Then
\[
g_0^*(N,k) = \tfrac{k-1}{12}N - \tfrac12 + c_2(k)\big( \tfrac{-1}N \big) + c_3(k)\big( \tfrac{-3}N \big),
\]
where $(\frac dN)$ is Kronecker's extension of the Legendre symbol. In particular, $g_0^*(N,k)$ depends on the residue class of $N$ modulo 12, but not on the prime factorization of $N$.
\label{gekeler.cor}
\end{cor}

We remark that the symbols $(\frac{-1}N)$ and $(\frac{-3}N)$ could also be represented by the nonprincipal characters $\chi_{-4}$ and $\chi_{-3}$ modulo 4 and 3, respectively. Gekeler used a proof by induction on the number of prime factors of $N$, which yielded a formula more complicated than, but equivalent to, the formula in Corollary \ref{gekeler.cor}. The corollary follows immediately from Theorem \ref{g0star.formula.theorem} by noting that $\delta(\frac k2)\delta(N)=0$ under the hypothesis $(k,N)\ne(2,1)$ and that $s_0^*(p) = \nu_\infty^*(p) = 1$, $\nu_2^*(p) = (\tfrac{-1}p)$, and $\nu_3^*(p) = (\tfrac{-3}p)$ for every prime $p$.

The situation is exactly the same for modular forms on $\Gamma_1(N)$: although the dimensions of spaces of cusp forms on $\Gamma_1(N)$ are well-understood, the dimensions of the corresponding spaces of newforms are more mysterious. Let $\Sk1N$ denote the space of cusp forms on $\Gamma_1(N)$ of weight $k$ and $\Skplus1N$ the space of newforms on $\Gamma_1(N)$ of weight $k$. Let $g_1(N,k)$ and $g_1^+(N,k)$ denote the dimensions of $\Sk1N$ and  $\Skplus1N$, respectively.  Also let $g_1^*(N,k)$ denote the number of nonisomorphic automorphic representations associated with $\Sk1N$. The method of proof of Theorems \ref{g0plus.formula.theorem} and \ref{g0star.formula.theorem} can also be used to establish formulas for $g_1^+(N,k)$ and $g_1^*(N,k)$ for any integer $k\ge2$ (not necessarily even). Since the expressions that result are slightly more complicated than the above expressions for $g_0^+(N,k)$ and $g_0^*(N,k)$, we defer the statements of the formulas to Theorems \ref{g1plus.formula.theorem} and \ref{g1star.formula.theorem} in Section \ref{Gamma1.section}. The complications arise because the most natural formula for $g_1(N,k)$ holds only for $N\ge5$; the presence of elliptic points and irregular cusps corresponding to $\Gamma_1(N)$ for $1\le N\le 4$ causes $g_1(N,k)$ to be somewhat different for these small values of $N$. Unfortunately, the behavior of $g_1^+(N,k)$ and $g_1^*(N,k)$ depends on the values of $g_1(N',k)$ for all divisors $N'$ of $N$, and so the exceptional cases $1\le N'\le4$ influence every single value of $g_1^+(N,k)$ and $g_1^*(N,k)$.

The explicit nature of the formulas in these theorems allows us to determine both the precise average orders and sharp asymptotic upper and lower bounds for these counting functions as well. The minimal and maximal orders of these functions are given in the next two theorems. Recall that $\gamma = \lim_{x\to\infty} \big( \sum_{n\le x} \frac1n - \log x \big) \approx 0.577216$ is Euler's constant.

%:Theorem - asymptotic bounds g0
\begin{theorem}
Uniformly for all even integers $k\ge2$ and all integers $N\ge1$, we have:
\begin{enumerate}
\item $\frac{k-1}{12}N + O(\sqrt N\log\log N) < g_0(N,k) < \frac{e^\gamma(k-1)}{2\pi^2}N\log\log N + O(N);$ %
\item $\frac{k-1}{2\pi^2}N + O\big( \frac{\phi(N)}{\sqrt N} \big) < g_0^*(N,k) < \frac{k-1}{12}N + O(1)$;
\item $\frac {A_0^+(k-1)}{12}\phi(N) + O(\sqrt N) < g_0^+(N,k) < \frac{k-1}{12}\phi(N) + O(2^{\omega(N)})$, where
\begin{equation}
A_0^+ = \prod_p \big( 1-\tfrac1{p^2-p} \big) \approx 0.373956.
\label{A0plus.def}
\end{equation}
Moreover, if $N$ is not a perfect square, then the lower bound can be improved to

\centerline{$\frac{A_0^+(k-1)}{12}\phi(N) + O(2^{\omega(N)}) < g_0^+(N,k)$.}
\end{enumerate}
\label{asymptotic.bounds.g0.theorem}
\end{theorem}

The product defining $A_0^+$ in equation \eqref{A0plus.def} is an infinite product over all prime numbers~$p$. The upper bounds in Theorem \ref{asymptotic.bounds.g0.theorem} imply in particular that both $g_0^*(N,k)$ and $g_0^+(N,k)$ are bounded above by a constant multiple of $kN$, in contrast to the size of $g_0(N,k)$ itself which can be as large as a constant multiple of $kN\log\log N$. Theorem \ref{asymptotic.bounds.g0.theorem} is stronger and more general than \cite[Proposition B.1]{HK}, which appeared as this manuscript was being prepared.
%We remark that we could easily derive a much better error term for the upper bound in part (a) of the theorem, such as $O(N\exp(-c\sqrt{\log\log N}))$ for some positive constant $c$, by appealing to the prime number theorem.

%:Theorem - asymptotic bounds g1
\begin{theorem}
Uniformly for all integers $k\ge2$ and all integers $N\ge1$, we have:
\begin{enumerate}
\item $\frac{k-1}{4\pi^2}N^2 + O(N\tau(N)+k) < g_1(N,k) < \frac{k-1}{24}N^2 + O(k)$;
\item $\frac{A_1^*(k-1)}{24}N^2 + O(N\tau(N)+k) < g_1^*(N,k) \le g_1(N,k)$, where
\begin{equation}
A_1^* = \prod_p \big( 1-\tfrac2{p^2} \big) \approx 0.322634;
\label{A1star.def}
\end{equation}
\item $\frac{A_1^+(k-1)}{24}N^2 + O(N\tau(N)+k) < g_1^+(N,k) \le g_1^*(N,k)$, where
\begin{equation}
A_1^+ = \prod_p \big( 1-\tfrac3{p^2} \big) \approx 0.125487.
\label{A1plus.def}
\end{equation}
\end{enumerate}
\label{asymptotic.bounds.g1.theorem}
\end{theorem}

To judge the quality of these error terms, recall that both $2^{\omega(N)}$ and $\tau(N)$ are $O(N^\ep)$ for any fixed $\ep>0$. Although Theorems \ref{asymptotic.bounds.g0.theorem}(a) and \ref{asymptotic.bounds.g1.theorem}(a) are easy consequences of the well-known formulas for $g_0(N,k)$ and $g_1(N,k)$, the bounds contained therein do not seem to have been recorded in the literature. We remark that all of the bounds given in Theorems \ref{asymptotic.bounds.g0.theorem} and \ref{asymptotic.bounds.g1.theorem} are best possible; the proofs of these theorems in Section~\ref{min.max.orders.section} are easily modified to produce sequences of values of $N$ asymptotically attaining the indicated upper and lower bounds.

We turn now to the question of the average orders of these various functions. Recall that a function $f(n)$ is said to have average order $g(n)$ if
\[\textstyle
\sum_{n\le x} f(n) \sim \sum_{n\le x} g(n),
\]
meaning that the quotient of the two sides approaches 1 as $x$ tends to infinity. It turns out that the average orders of the counting functions associated with $\Gamma_0(N)$ are explicit constant multiples of $N$.

%:Theorem - average orders g_0
\begin{theorem}
Fix an even integer $k\ge2$.
\begin{enumerate}
\item The average order of $g_0(N,k)$ is $5(k-1)N/4\pi^2$.
\item The average order of $g_0^*(N,k)$ is $15(k-1)N/2\pi^4$.
\item The average order of $g_0^+(N,k)$ is $45(k-1)N/\pi^6$.
\end{enumerate}
\label{average.order.g0.theorem}
\end{theorem}

The average orders of the counting functions associated with $\Gamma_1(N)$ depend on the special value $\zeta(3) = \sum_{n=1}^\infty n^{-3}$ of the Riemann zeta-function.

%:Theorem - average orders g_1
\begin{theorem}
Fix an integer $k\ge2$.
\begin{enumerate}
\item The average order of $g_1(N,k)$ is $(k-1)N^2/24\zeta(3)$.
\item The average order of $g_1^*(N,k)$ is $(k-1)N^2/24\zeta(3)^2$.
\item The average order of $g_1^+(N,k)$ is $(k-1)N^2/24\zeta(3)^3$.
\end{enumerate}
\label{average.order.g1.theorem}
\end{theorem}

Another natural quantity to consider is the relative number of newforms with the spaces of cusp forms on $\Sk0N$ and $\Sk1N$. To measure this proportion, define 
\[
\rho_0(N,k) = \begin{cases}
g_0^+(N,k)/g_0(N,k), &\text{if } g_0(N,k)>0, \\
1, &\text{if } g_0(N,k)=0 \\
\end{cases}
\]
and similarly for $\rho_1(N,k)$. We are able to establish asymptotically sharp lower bounds for $\rho_0(N,k)$ and $\rho_1(N,k)$, as well as determine their average orders.

%:Theorem - bounds for the rhos
\begin{theorem}
Uniformly for all integers $k\ge2$ and all integers $N\ge1$, we have:
\begin{enumerate}
\item $\frac{A_0^+\pi^2}{6e^{2\gamma}(\log\log N)^2} + O\big( \frac1{(\log\log N)^3} \big) < \rho_0(N,k) \le 1$, where $A_0^+$ is defined in equation~\eqref{A0plus.def};
\item $\frac{A_1^+\pi^2}6 + O\big( \frac1{\log N\log\log N} + \frac kN \big) < \rho_1(N,k) \le 1$, where $A_1^+$ is defined in equation~\eqref{A1plus.def}.
\end{enumerate}
\label{asymptotic.bounds.rhos.theorem}
\end{theorem}

Note that $\frac{A_1^+\pi^2}6 \approx 0.206418$; we deduce from the lower bound in Theorem \ref{asymptotic.bounds.rhos.theorem}(b) that when $N$ is large enough with respect to $k$, it always the case that at least 20\% of the weight-$k$ cusp forms on $\Gamma_1(N)$ are newforms.
% EXCEPTIONS: whenever $g_1(N,k)=0$; and $N=2$, $k=12$ or 24; and $N=4$ or 6, $k$ even (except for a few small $k$). Explicit constants and a computation could show that these are the only exceptions.

%:Theorem - average orders for the rhos
\begin{theorem}
Fix an integer $k\ge2$.
\begin{enumerate}
\item If $k$ is even, then the average order of $\rho_0(N,k)$ is
\begin{equation}
B_0 = \prod_p \big( 1-\tfrac1p \big) \big( 1+\tfrac1p \big)^{-1} \big( 1+\tfrac2p - \tfrac1{p^4} -\tfrac1{p^5} \big) \approx 0.444301.
\label{B0.def}
\end{equation}
\item The average order of $\rho_1(N,k)$ is
\begin{equation}
B_1 = \prod_p \big( 1+\tfrac1p \big)^{-1} \big( 1 + \tfrac1p - \tfrac2{p^3} - \tfrac2{p^4} - \tfrac2{p^5} + \tfrac1{p^6} + \tfrac1{p^7} + \tfrac1{p^8} \big) \approx 0.652036.
\label{B1.def}
\end{equation}
\end{enumerate}
\label{average.order.rhos.theorem}
\end{theorem}

In Section \ref{main.theorems.proofs.section}, we prove the main formulas for $g_0^+(N,k)$ and $g_0^*(N,k)$ given in Theorems \ref{g0plus.formula.theorem} and \ref{g0star.formula.theorem}. Subsequently, we investigate the analogous functions for modular forms on $\Gamma_1(N)$ in Section \ref{Gamma1.section}, culminating in the statements and proofs of Theorems \ref{g1plus.formula.theorem} and \ref{g1star.formula.theorem}. Sections \ref{explicit.section} and \ref{up.to.100.section} are devoted to the explicit inequalities in Theorems \ref{bennett.bound.theorem} and \ref{bennett.inequality.theorem} and to computational resuts concerning $g_0^+(N,2)$ and $g_0(N,2)$. We finish by establishing the asymptotic inequalities of Theorems \ref{asymptotic.bounds.g0.theorem}, \ref{asymptotic.bounds.g1.theorem}, and \ref{asymptotic.bounds.rhos.theorem} in Section \ref{min.max.orders.section} and the average-order results of Theorems \ref{average.order.g0.theorem}, \ref{average.order.g1.theorem}, and \ref{average.order.rhos.theorem} in Section \ref{average.section}.

\section{Notation and proof of Theorems
\ref{g0plus.formula.theorem} and \ref{g0star.formula.theorem}}
\label{main.theorems.proofs.section}

The dimensions of the spaces of weight-$k$ cusp forms on $\Gamma_0(N)$ are well-known for positive even integers $k$. The following proposition gives a formula for these dimensions, phrased in the way that is most convenient for our purposes.

%:Proposition - formula for g_0
\begin{prop}
For any even integer $k\ge2$ and any integer $N\ge1$, we have
\[
g_0(N,k) = \tfrac{k-1}{12}Ns_0(N) - \tfrac12\nu_\infty(N)+c_2(k) \nu_2(N) +c_3(k) \nu_3(N) + \delta\big( \tfrac k2 \big),
\]
where the functions $s_0$, $\nu_\infty$, $\nu_2$, $\nu_3$, $c_2$, and $c_3$ are defined in Definitions \ref{s0def}--\ref{nu3def} below and Definitions \ref{c2def}--\ref{c3def} above.
\label{g0.formula.prop}
\end{prop}

%:(Definitions needed to state the g_0 formula)
The definitions of the four new functions in the statement of Proposition \ref{g0.formula.prop} are as follows.

\begin{definition}
$s_0$ is the multiplicative function satisfying $s_0(p^\alpha) = 1+\tfrac1p$\, for all $\alpha\ge1$.
\label{s0def}
\end{definition}

\begin{definition}
$\nu_\infty$ is the multiplicative function satisfying

\centerline{$\displaystyle
\nu_\infty(p^\alpha) = \begin{cases}
2p^{(\alpha-1)/2}, &\text{if $\alpha$ is odd}, \\
p^{\alpha/2} + p^{\alpha/2-1}, &\text{if $\alpha$ is even}.
\end{cases}
$}
\label{nuinftydef}
\end{definition}

\begin{definition}
$\nu_2$ is the multiplicative function satisfying:
\begin{itemize}
\item $\nu_2(2)=1$\, and $\nu_2(2^\alpha)=0$ for $\alpha\ge2$;
\item if $p\equiv1\mod4$ then $\nu_2(p^\alpha)=2$ for $\alpha\ge1$;
\item if $p\equiv3\mod4$ then $\nu_2(p^\alpha)=0$ for $\alpha\ge1$.
\end{itemize}
\label{nu2def}
\end{definition}

\begin{definition}
$\nu_3$ is the multiplicative function satisfying:
\begin{itemize}
\item $\nu_3(3)=1$\, and $\nu_3(3^\alpha)=0$ for $\alpha\ge2$;
\item if $p\equiv1\mod3$ then $\nu_3(p^\alpha)=2$ for $\alpha\ge1$;
\item if $p\equiv2\mod3$ then $\nu_3(p^\alpha)=0$ for $\alpha\ge1$.
\end{itemize}
\label{nu3def}
\end{definition}

\begin{proof}[Proof of Proposition \ref{g0.formula.prop}]
The facts invoked in this proof can be found in many sources; we follow the exposition in Miyake \cite{miyake}. For now we assume that $N\ge2$. We begin by remarking that the multiplicative function $\nu_\infty(N)$ denotes the number of (inequivalent) cusps of $\Gamma_0(N)$ and that the multiplicative functions $\nu_j(N)$ denote the numbers of (inequivalent) elliptic points of $\Gamma_0(N)$ of order~$j$. Formulas for these numbers are given in \cite[Theorem 4.2.7]{miyake} in the form
\begin{equation}
\nu_\infty(N) = \sum_{d\mid n} \phi\big( (d,\tfrac nd) \big) = \prod_{p^\alpha\exdiv N} \bigg\{ \sum_{\beta=0}^\alpha \phi\big( p^{\min\{\beta,\alpha-\beta\}} \big) \bigg\}
\label{miyake.nuinfty}
\end{equation}
and
\begin{equation}
\nu_2(N) = \begin{cases}
0, &\text{if }4\mid n, \\
\prod\limits_{p\mid n} \big( 1 + (\tfrac{-1}p) \big), &\text{otherwise;}
\end{cases} \quad\ 
\nu_3(N) = \begin{cases}
0, &\text{if }9\mid n, \\
\prod\limits_{p\mid n} \big( 1 + (\tfrac{-3}p) \big), &\text{otherwise.}
\end{cases}
\label{miyake.nu2.nu3}
\end{equation}
Here again the symbol $(\frac ap)$ is Kronecker's extension of the Legendre symbol. It is easily verified that the formulas for $\nu_2$ and $\nu_3$ in equation \eqref{miyake.nu2.nu3} are equivalent to the formulas in Definitions \ref{nu2def} and \ref{nu3def}. It is also easily verified that since $\alpha\ge1$,
\[
\sum_{\beta=0}^\alpha \phi\big( p^{\min\{\beta,\alpha-\beta\}} \big) = 2 + (p-1) \sum_{\beta=1}^{\alpha-1} p^{\min\{\beta,\alpha-\beta\}-1} = \begin{cases}
2p^{(\alpha-1)/2}, &\text{if $\alpha$ is odd}, \\
p^{\alpha/2} + p^{\alpha/2-1}, &\text{if $\alpha$ is even},
\end{cases}
\]
and so the formula in equation \eqref{miyake.nuinfty} is the same as the formula in Definition \ref{nuinftydef}.

Next, if we let $g_N$ denote the genus of the (compactified) quotient of the upper half-plane by $\Gamma_0(N)$, then we have the formula \cite[Theorem 4.2.11]{miyake}
\begin{equation}
g_N = \tfrac{\mu_N}{12} - \tfrac{\nu_\infty(N)}2 - \tfrac{\nu_2(N)}4 -\tfrac{\nu_3(N)}3+1,
\label{gN.formula}
\end{equation}
where $\mu_N$ is the index of $\overline\Gamma_0(N)$ in $\overline{SL}_2(\Z)$, and  $\overline G$ denotes the quotient of the group $G$ by its center. According to \cite[Theorem 4.2.5]{miyake},
\[
\mu_N = N\prod_{p\mid N} \big( 1+\tfrac1p \big) = Ns_0(N)
\]
as defined in Definition \ref{s0def}.

Now the dimension $g_0(N,k)$ of the space of weight-$k$ cusp forms on $\Gamma_0(N)$ can be calculated from this information by the Riemann--Roch theorem. From \cite[Theorem 2.5.2]{miyake} we see that $g_0(N,2)=g_N$ and
\[
g_0(N,k) = (k-1)(g_N-1) + \big( \tfrac k2-1 \big) \nu_\infty(N) + \sum_{j\ge2} \big\lfloor \tfrac k2\big( 1-\tfrac1j \big) \big\rfloor \nu_j(N)
\]
for every even integer $k\ge4$. Only the terms $j=2,3$ are present in the sum due to \cite[Lemma 4.2.6]{miyake}, and so the equation for $g_0(N,k)$ becomes
\[
g_0(N,k) = (k-1)(g_N-1) + \big( \tfrac k2-1 \big) \nu_\infty(N) + \big\lfloor \tfrac k4 \big\rfloor \nu_2(N) + \big\lfloor \tfrac k3 \big\rfloor \nu_3(N).
\]
Combining this with the formula \eqref{gN.formula} and collecting the multiples of $\nu_\infty(N)$, $\nu_2(N)$, and $\nu_3(N)$ yields
\begin{equation}
g_0(N,k) = \tfrac{k-1}{12}Ns_0(N) - \tfrac12 \nu_\infty(N) + \big( \tfrac14 - \tfrac k4 + \big\lfloor \tfrac k4 \big\rfloor \big) \nu_2(N) + \big( \tfrac13 - \tfrac k3 + \big\lfloor \tfrac k3 \big\rfloor \big) \nu_3(N),
\label{g0Nk.with.cks.written.out}
\end{equation}
which is the same as the assertion of the proposition (when $k\ge4$) in light of the definitions \ref{c2def} and \ref{c3def} of $c_2$ and $c_3$. It is easily checked that the formula holds for $k=2$ as well. Finally, all of this discussion assumed that $N\ge2$, but the special case $N=1$ is worked through in detail in \cite[Section 4.1]{miyake}, and the formula \cite[Corollary 4.1.4]{miyake} can be seen to agree with the assertion of the proposition as well.
\end{proof}

%From \cite{miyake} or \cite{Shimura} we know that for even integers $k\ge2$,
%\begin{equation}
%g_0(N,k) = \tfrac{k-1}{12}Ns_0(N) - \tfrac12\nu_\infty(N) +c_2(k) \nu_2(N) +c_3(k) \nu_3(N) + \delta\big( \tfrac k2 \big).
%\label{g0nk.formula}
%\end{equation}
%In other words, when $k\ge4$ is even,
%\begin{multline*}
%g_0(N,k) = \tfrac{k-1}{12}Ns_0(N) - \tfrac12\nu_\infty(N) \\
%+ \begin{cases}
%\hskip3.3mm\frac14\nu_2(N), &\text{if } k\equiv2\mod4 \\
%-\frac14\nu_2(N), &\text{if } k\equiv0\mod4 \\
%\end{cases} \!\Bigg\}
%+ \begin{cases}
%\hskip3.3mm\frac13\nu_3(N), &\text{if } k\equiv2\mod6, \\
%\hskip13mm0, &\text{if } k\equiv4\mod6, \\
%-\frac13\nu_3(N), &\text{if } k\equiv0\mod6. \\
%\end{cases}
%\end{multline*}

We may now prove Theorems \ref{g0plus.formula.theorem} and \ref{g0star.formula.theorem}.

%:Proof of Theorem \ref{g0plus.formula.theorem}
\begin{proof}[Proof of Theorem \ref{g0plus.formula.theorem}]
If $f(z)$ is a cusp form on $\Gamma_0(d)$, then $f(mz)$ is a cusp form on $\Gamma_0(N)$ for any multiple $N$ of $dm$. Therefore for every triple $(m,d,N)$ of positive integers such that $dm\mid N$, we have an injection $i_{m,d,N}: \Sk0d \to \Sk0N$ defined by $i_{m,d,N}(f)(z) = f(mz)$. As shown by Atkin and Lehner \cite{AL}, we may write
\begin{equation}
\Sk0N = \bigoplus_{d\mid N} \bigoplus_{m\mid N/d} i_{m,d,N}\big(\Skplus0d\big)
\label{atkin.lehner.decomposition}
\end{equation}
(in fact, summands corresponding to distinct divisors $d$ are orthogonal with respect to the Petersson inner product). In particular, the dimensions of these spaces satisfy
\begin{equation}
g_0(N,k) = \sum_{d\mid N} \sum_{m\mid N/d} g_0^+(d,k) = \sum_{d\mid N} g_0^+(d,k) \tau(N/d).
\label{g0.in.terms.of.g0plus}
\end{equation}

This equation can be written more simply using the Dirichlet convolution
\begin{equation}
f*g(n) = \sum_{d\mid n} f(d) g(n/d).
\label{Dirichlet.convolution}
\end{equation}
Recall that the set of arithmetic functions $f\colon\N\to\C$ forms a ring under the usual addition of functions and the Dirichlet convolution as the multiplication operation, with the function $\delta$ defined in equation \eqref{delta.def} as the multiplicative identity. In fact, the set of multiplicative functions forms a multiplicative subgroup---the Dirichlet convolution of two multiplicative functions $f,g$ is again multiplicative. Indeed, the values of $f*g$ on prime powers can be computed easily from the values of $f$ and $g$ on prime powers using the identity
\begin{equation}
f*g(p^\alpha) = \sum_{\beta=0}^\alpha f(p^\beta) g(p^{\alpha-\beta}),
\label{primepower.convolution}
\end{equation}
which is a special case of equation \eqref{Dirichlet.convolution}. We also remark that the characteristic property of the M\"obius $\mu$ function, often phrased as the M\"obius inversion formula, is that it is the inverse (under Dirichlet convolution) of the function $1(n)$ that takes the value 1 at all positive integers:
\[
(\mu*1)(n) = \sum_{d\mid n} \mu(d) = \delta(n).
\]

Now in this notation, equation \eqref{g0.in.terms.of.g0plus} says simply that $g_0 = g_0^+ * \tau$ for every fixed $k$. Define $\lambda$ to be the inverse (under Dirichlet convolution) of $\tau$. Since $\tau=1*1$, we see that $\lambda = \mu*\mu$. Equivalently, $\lambda$ is the multiplicative function satisfying
\begin{equation}
\lambda(p)={-2}, \quad \lambda(p^2)=1, \quad \lambda(p^\alpha)=0 \text{ for }\alpha\ge3,
\label{lambda.def}
\end{equation}
as can be seen by applying the formula \eqref{primepower.convolution} with $f=g=\mu$. It follows that $g_0^+ = g_0 * \lambda$ for every fixed $k$, that is,
\begin{equation}
g_0^+(N,k) = \sum_{d\mid N} g_0(d,k) \lambda(N/d).
\label{g0plus.in.terms.of.g0}
\end{equation}
However, since $g_0^+(N,k)$ is a linear combination of multiplicative functions of $N$ (with coefficients depending on $k$), it is more natural to take the convolution of $\lambda$ with the right-hand side of the formula given in Proposition \ref{g0.formula.prop}. We obtain
\begin{multline*}
g_0^+(N,k) = \tfrac{k-1}{12}Ns_0(N)*\lambda(N) - \tfrac12 (\nu_\infty*\lambda)(N) \\
+ c_2(k) (\nu_2*\lambda)(N) + c_3(k) (\nu_3*\lambda)(N) + \delta\big( \tfrac k2 \big) (1*\lambda)(N).
\end{multline*}
We immediately note that $1*\lambda=1*\mu*\mu=\mu$. Furthermore, the functions $\nu_\infty*\lambda$, $\nu_2*\lambda$, and $\nu_3*\lambda$ are all multiplicative; by using the formula \eqref{primepower.convolution} we see that they are equal to the functions $\nu_\infty^+$, $\nu_2^+$, and $\nu_3^+$ defined in Definitions \ref{nuinftyplusdef}--\ref{nu3plusdef}. Finally, it can be verified using \eqref{primepower.convolution} that
\[
p^\alpha s_0(p^\alpha)*\lambda(p^\alpha) = \sum_{\beta=0}^\alpha p^\beta s_0(p^\beta) \lambda(p^{\alpha-\beta}) = p^\alpha s_0^+(p^\alpha),
\]
where $s_0^+$ is defined in Definition \ref{s0plusdef}; therefore the multiplicative function $Ns_0(N)*\lambda(N)$ is equal to $Ns_0^+(N)$. This establishes the theorem.
\end{proof}

%:Proof of Theorem \ref{g0star.formula.theorem}
\begin{proof}[Proof of Theorem \ref{g0star.formula.theorem}]
The spaces of cusp forms $\Sk0N$ have bases consisting of modular forms that are eigenforms for all but finitely many Hecke operators. An isomorphism class of automorphic representations corresponds to an equivalence class of eigenforms, where two eigenforms are equivalent if all but finitely many Hecke operators act upon them with the same eigenvalues, or equivalently if both eigenforms are the image of the same newform under two injections $i_{m_1,d,N}$ and  $i_{m_2,d,N}$. Therefore, if we define the subspace $\Skstar0N$ of $\Sk0N$ to be
\begin{equation}
\Skstar0N = \bigoplus_{d\mid N} i_{1,d,N}\big(\Skplus0d\big),
\label{Skstar0N.definition}
\end{equation}
then the dimension of $\Skstar0N$ can be interpreted as the number of nonisomorphic automorphic representations associated with $\Sk0N$, which we have denoted by $g_0^*(N,k)$. From here, the proof is very similar to the proof of Theorem \ref{g0plus.formula.theorem}. The dimensions of these spaces satisfy
\begin{equation*}
g_0^*(N,k) = \sum_{d\mid N} g_0^+(d,k);
%\label{g0star.in.terms.of.g0plus}
\end{equation*}
in other words, $g_0^*$ is simply the convolution $g_0^+*1$ for every fixed $k$. We saw in the proof of Theorem \ref{g0plus.formula.theorem} that $g_0^+ = g_0*\lambda$ for every fixed $k$, and hence $g_0^*= g_0*\lambda*1 = g_0*\mu$, that~is,
\begin{equation*}
g_0^*(N,k) = \sum_{d\mid N} g_0(d,k) \mu(N/d).
%\label{g0star.in.terms.of.g0}
\end{equation*}
Again, since $g_0^+(N,k)$ is a linear combination of multiplicative functions of $N$ (with coefficients depending on $k$), it is natural to use Proposition \ref{g0.formula.prop} to write
\begin{multline*}
g_0^*(N,k) = \tfrac{k-1}{12}Ns_0(N)*\mu(N) - \tfrac12 (\nu_\infty*\mu)(N) \\
+ c_2(k) (\nu_2*\mu)(N) + c_3(k) (\nu_3*\mu)(N) + \delta\big( \tfrac k2 \big) (1*\mu)(N).
\end{multline*}
We immediately note that $1*\mu=\delta$. Furthermore, the functions $\nu_\infty*\mu$, $\nu_2*\mu$, and $\nu_3*\mu$ are all multiplicative; by using the formula \eqref{primepower.convolution} we see that they are equal to the functions $\nu_\infty^+$, $\nu_2^+$, and $\nu_3^+$ defined in Definitions \ref{nuinftystardef}--\ref{nu3stardef}. Finally, using \eqref{primepower.convolution} we verify that
\[
p^\alpha s_0(p^\alpha)*\mu(p^\alpha) = \sum_{\beta=0}^\alpha p^\beta s_0(p^\beta) \mu(p^{\alpha-\beta}) = p^\alpha s_0(p^\alpha) - p^{\alpha-1} s_0(p^{\alpha-1}) = p^\alpha s_0^*(p^\alpha),
\]
where $s_0^*$ is defined in Definition \ref{s0stardef}; therefore the multiplicative function $Ns_0(N)*\mu(N)$ is equal to $Ns_0^*(N)$. This establishes the theorem.
\end{proof}

\section{formulas for $g_1^+$ and $g_1^*$}
\label{Gamma1.section}

In this section we state and prove formulas for modular forms on $\Gamma_1(N)$ that are analogous to Theorems \ref{g0plus.formula.theorem} and \ref{g0star.formula.theorem}.

%:Theorem - formula for g_1^+
\begin{theorem}
For any integer $k\ge2$ and any integer $N\ge1$, we have
\[
g_1^+(N,k) = \tfrac{k-1}{24}N^2s_1^+(N) - \tfrac14u^+(N) + \delta\big( \tfrac k2 \big)\mu(N) + \sum_{\substack{1\le i\le 4 \\ i\mid N}} b_i(k) \lambda(N/i),
\]
where the functions $s_1^+$, $u^+$, $b_1$, $b_2$, $b_3$, and $b_4$ are defined in Definitions \ref{s1plusdef}--\ref{bdefs} below.
\label{g1plus.formula.theorem}
\end{theorem}

%:(Definitions needed to state the g_1^+ formula)
Recall that the multiplicative function $\lambda=\mu*\mu$ was defined in equation \eqref{lambda.def} above. The definitions of the six functions in the statement of Theorem \ref{g1plus.formula.theorem} are as follows.

\begin{definition}
$s_1^+$ is the multiplicative function satisfying

\centerline{$
s_1^+(p) = 1-\tfrac3{p^2},\;
s_1^+(p^2) = 1-\tfrac3{p^2}+\tfrac3{p^4},\;
\text{and } s_1^+(p^\alpha) = \big(1-\tfrac1{p^2}\big)^3 \text{ for } \alpha\ge3.
$}
\label{s1plusdef}
\end{definition}

\begin{definition}
$u^+$ is the multiplicative function satisfying $u^+(p)=2p-4$,\,  $u^+(p^2)=3p^2-8p+6$,\, and

\centerline{$
u^+(p^\alpha) = p^{\alpha-4}(p-1)^3((\alpha+1)p-\alpha+3) \text{ for } \alpha\ge3.
$}
\label{uplusdef}
\end{definition}

\begin{definition}
The functions $b_i(k)$ are defined as follows:
\begin{itemize}
\item $\displaystyle b_1(k) = \tfrac{(-1)^k(k-7)}{24} + \begin{cases}
c_2(k)+c_3(k), &\text{if $k$ is even,} \\ 0, &\text{if $k$ is odd;} \\
\end{cases}$
\item $b_2(k) = \frac12 \big( (-1)^k \big\lfloor \frac k4-1 \big\rfloor + c_2(k) \big)$;
\item $b_3(k) = c_3(k)$;
\item $b_4(k) = -c_2(2k)$.
\end{itemize}
\label{bdefs}
\end{definition}

There are many equivalent ways to write the formulas defining the functions $b_i(k)$. Our choices were motivated by the desire to make the sizes of the functions $b_i(k)$ as $k$ grows immediately apparent, knowing that the functions $c_2(k)$ and $c_3(k)$ are bounded in absolute value by $\frac12$.

%:Theorem - formula for g_1^*
\begin{theorem}
For any integer $k\ge2$ and any integer $N\ge1$, we have
\[
g_1^*(N,k) = \tfrac{k-1}{24}N^2s_1^*(N) - \tfrac14u^*(N) + \delta\big( \tfrac k2 \big)\delta(N) + \sum_{\substack{1\le i\le 4 \\ i\mid N}} b_i(k) \mu(N/i),
\]
where the functions $s_1^*$, $u^*$, $b_1$, $b_2$, $b_3$, and $b_4$ are defined in Definitions \ref{s1stardef}--\ref{ustardef} below and Definition \ref{bdefs} above.
\label{g1star.formula.theorem}
\end{theorem}

%:(Definitions needed to state the g_1^* formula)
The definitions of the two new functions in the statement of Theorem \ref{g1star.formula.theorem} are as follows.

\begin{definition}
$s_1^*$ is the multiplicative function satisfying

\centerline{$
s_1^*(p) = 1-\tfrac2{p^2}\;
\text{and } s_1^*(p^\alpha) = \big(1-\tfrac1{p^2}\big)^2 \text{ for } \alpha\ge2.
$}
\label{s1stardef}
\end{definition}

\begin{definition}
$u^*$ is the multiplicative function satisfying $u^*(p)=2p-4$\, and $u^*(p^\alpha) = p^{\alpha-3}(p-1)^2((\alpha+1)p-\alpha+2)$ for $\alpha\ge2$.
\label{ustardef}
\end{definition}

As in the previous section, our starting point is a formula for $g_1(n,k)$, the dimension of the space $\Sk1N$ of weight-$k$ modular forms on $\Gamma_1(N)$.

%:Proposition - formula for g_1
\begin{prop}
For any integer $k\ge2$ and any integer $N\ge1$, we have
\begin{equation}
g_1(N,k) = \tfrac{k-1}{24}N^2s_1(N) - \tfrac14u(N) + \delta\big( \tfrac k2 \big) + \sum_{\substack{1\le i\le 4 \\ i\mid N}} b_i(k) \delta(N/i),
\label{g1.formula.eq}
\end{equation}
where the functions $s_1$, $u$, $b_1$, $b_2$, $b_3$, and $b_4$ are defined in Definitions \ref{s1def}--\ref{udef} below and Definition \ref{bdefs} above.
\label{g1.formula.prop}
\end{prop}

%:(Definitions needed to state the g_1 formula)
The definitions of the two new functions in the statement of Proposition \ref{g1.formula.prop} are as follows.

\begin{definition}
$s_1$ is the multiplicative function satisfying $s_1(p^\alpha) = 1-\tfrac1{p^2}$\, for all $\alpha\ge1$.
\label{s1def}
\end{definition}

\begin{definition}
$u$ is the multiplicative function satisfying

\centerline{
$u(p^\alpha) = p^{\alpha-2}(p-1)((\alpha+1)p-\alpha+1)$\, for all $\alpha\ge1$.}
\label{udef}
\end{definition}

\begin{proof}
As in the proof of Proposition \ref{g0.formula.prop}, our main task is simply to gather together the known facts about $\Gamma_1(N)$. For now we assume that $N\ge5$. In this case, by \cite[Theorem 4.2.9]{miyake}, we know both that $\Gamma_1(N)$ has no elliptic elements and that the number of (inequivalent) cusps of $\Gamma_1(N)$ is given by the formula $\tfrac12 \sum_{d\mid n} \phi(d)\phi(n/d)$. We calculate that
\begin{align*}
\sum_{d\mid n} \phi(d)\phi(n/d) &= \prod_{p^\alpha\exdiv n} \sum_{d\mid p^\alpha} \phi(d)\phi(p^\alpha/d) \\
&= \prod_{p^\alpha\exdiv n} \sum_{\beta=0}^\alpha \phi(p^\beta)\phi(p^{\alpha-\beta}) \\
&= \prod_{p^\alpha\exdiv n} \big( 2p^{\alpha-1}(p-1) + (\alpha-1)p^{\alpha-2}(p-1)^2 \big) \\
&= \prod_{p^\alpha\exdiv n} p^{\alpha-2}(p-1)((\alpha+1)p-\alpha+1).
\end{align*}
Thus this expression for the number of cusps is nothing other than $\frac12u(n)$ as defined in Definition \ref{udef}.

We now let $g_N$ denote the genus of the quotient of the upper half-plane by $\Gamma_1(N)$ and $\mu_N$ the index of $\overline\Gamma_1(N)$ in $\overline{SL}_2(\Z)$, superceding the notation in the proof of Proposition \ref{g0.formula.prop}. From \cite[Theorem 4.2.5]{miyake}, we have that
\[
\mu_N = \tfrac{\phi(N)}2 \cdot N\prod_{p\mid N}\big( 1+\tfrac1p \big) = \tfrac12 N^2 \prod_{p\mid N} \big( 1-\tfrac1{p^2} \big) = \tfrac12 N^2 s_1(N)
\]
according to Definition \ref{s1def}. The formula \eqref{gN.formula} then becomes
\[
g_N = \tfrac{N^2 s_1(N)}{24} - \tfrac{u(N)}4 + 1.
\]
Using \cite[Theorem 2.5.2]{miyake} again, we discover that $g_1(N,2) = g_N$ and that for even $k\ge4$,
\[
g_1(N,k) = \tfrac{k-1}{24}N^2s_1(N) - \tfrac14u(N)
\]
in analogy with equation \eqref{g0Nk.with.cks.written.out}. We may combine these two facts into the single equation
\begin{equation}
g_1(N,k) = \tfrac{k-1}{24}N^2s_1(N) - \tfrac14u(N) + \delta\big(\tfrac k2\big),
\label{into.the.single.equation}
\end{equation}
in agreement with the assertion of the proposition (note that the sum in equation \eqref{g1.formula.eq} is zero when $N\ge5$). An appeal to \cite[Theorem 2.5.3]{miyake} shows that this equation holds when $k\ge3$ is odd as well. This establishes the proposition when $N\ge5$.

Unfortunately, the groups $\Gamma_1(N)$ for $1\le N\le4$ are exceptional, and the general formula just derived does not give the correct answer. When $1\le N\le4$ we have $\overline\Gamma_1(N) \cong \overline\Gamma_0(N)$, and so the true values of $g_1(N,k)$ for these small $N$ are equal to the values $g_1(N,k)$ when $k\ge2$ is even. Calculating these values explicitly from Proposition \ref{g0.formula.prop}, we have
\begin{align*}
g_1(1,k) &= \lfloor \tfrac k4 \rfloor + \lfloor \tfrac k3 \rfloor - \tfrac k2 + \delta(\tfrac k2) \\
g_1(2,k) &= \lfloor \tfrac k4 \rfloor - 1 + \delta(\tfrac k2) \\
g_1(3,k) &= \lfloor \tfrac k3 \rfloor - 1 + \delta(\tfrac k2) \\
g_1(4,k) &= \lfloor \tfrac {k-3}2 \rfloor - 1 + \delta(\tfrac k2)
\end{align*}
for even integers $k\ge2$. When $k\ge3$ is odd, we know that $g_1(1,k) = g_1(2,k)= 0$ since $\Gamma_1(1)=SL_2(\Z)$ and $\Gamma_1(2)$ both contain the matrix $\big( \genfrac{}{}{0pt}{}{-1\,}{\phantom{-}0\,} \genfrac{}{}{0pt}{}{\phantom{-}0}{-1} \big)$. By carefully working through the details in \cite[Section 4.2]{miyake}, we see that the above formulas for $g_1(3,k)$ and $g_1(4,k)$ are also correct when $k\ge3$ is odd. In other words, the formulas
\begin{align*}
g_1(1,k) &= \big( \tfrac{1+(-1)^k}2 \big) \big( \lfloor \tfrac k4 \rfloor + \lfloor \tfrac k3 \rfloor - \tfrac k2 \big) + \delta(\tfrac k2) \\
g_1(2,k) &= \big( \tfrac{1+(-1)^k}2 \big) \big( \lfloor \tfrac k4 \rfloor - 1 \big) + \delta(\tfrac k2) \\
g_1(3,k) &= \lfloor \tfrac k3 \rfloor - 1 + \delta(\tfrac k2) \\
g_1(4,k) &= \lfloor \tfrac {k-3}2 \rfloor + \delta(\tfrac k2)
\end{align*}
are valid for all $k\ge2$.

The formula \eqref{into.the.single.equation} gives the false values $\frac{k-7}{24} + \delta(\frac k2)$, $\frac{k-5}8 + \delta(\frac k2)$, $\frac{k-4}3 + \delta(\frac k2)$, and $\frac{2k-7}4 + \delta(\frac k2)$ for $g_1(1,k)$, $g_1(2,k)$, $g_1(3,k)$, and $g_1(4,k)$, respectively. One can check that
\begin{align*}
\big( \tfrac{1+(-1)^k}2 \big) \big( \lfloor \tfrac k4 \rfloor + \lfloor \tfrac k3 \rfloor - \tfrac k2 \big) - \tfrac{k-7}{24} &= b_1(k) \\
\big( \tfrac{1+(-1)^k}2 \big) \big( \lfloor \tfrac k4 \rfloor - 1 \big) - \tfrac{k-5}8 &= b_2(k) \\
\lfloor \tfrac k3 \rfloor - 1 - \tfrac{k-4}3 &= b_3(k) \\
\lfloor \tfrac {k-3}2 \rfloor - \tfrac{2k-7}4 &= b_4(k)
\end{align*}
using the definition \ref{bdefs} of the functions $b_i(k)$. Therefore we can write
\[
g_1(N,k) = \tfrac{k-1}{24}N^2s_1(N) - \tfrac14u(N) + \delta\big(\tfrac k2\big) + \begin{cases}
b_1(k), &\text{if } N=1, \\
b_2(k), &\text{if } N=2, \\
b_3(k), &\text{if } N=3, \\
b_4(k), &\text{if } N=4, \\
0, &\text{if } N\ge5, \\
\end{cases}
\]
which is equivalent to the assertion of the proposition for all $N\ge1$ and $k\ge2$.
\end{proof}

We may now prove Theorems \ref{g1plus.formula.theorem} and \ref{g1star.formula.theorem}.

%:Proof of formulae for g_1^+ and g_1^*
\begin{proof}[Proof of Theorems \ref{g1plus.formula.theorem} and \ref{g1star.formula.theorem}]
We proceed as in the proofs of Theorems \ref{g0plus.formula.theorem} and \ref{g0star.formula.theorem}. Again we have the Atkin--Lehner decomposition
\begin{equation*}
\Sk1N = \bigoplus_{d\mid N} \bigoplus_{m\mid N/d} i_{m,d,N}\big(\Skplus1n\big).
\end{equation*}
Calculating the dimensions of both sides yields
\begin{equation*}
g_1(N,k) = \sum_{d\mid N} \sum_{m\mid N/d} g_1^+(d,k) = \sum_{d\mid N} g_1^+(d,k) \tau(N/d).
\end{equation*}
This implies that $g_1^+ = g_1 * \lambda$ for every fixed $k$ (recall the definition \eqref{lambda.def} of the multiplicative function $\lambda$), that is,
\begin{equation*}
g_1^+(N,k) = \sum_{d\mid N} g_1(d,k) \lambda(N/d).
\end{equation*}
Using the formula for $g_1(N,k)$ given in Proposition \ref{g1.formula.prop}, this becomes
\[
g_1(N,k) = \tfrac{k-1}{24}N^2s_1(N)*\lambda(N) - \tfrac14(u*\lambda)(N) + \delta\big( \tfrac k2 \big) + \sum_{\substack{1\le i\le 4 \\ i\mid N}} b_i(k) \big( \delta(N/i)*\lambda(N) \big).
\]
We immediately note that the expression $\delta(N/i)*\lambda(N)$ equals simply $\lambda(N/i)$ in the case where $i$ divides $N$. Furthermore, the function $u*\lambda$ is multiplicative; by using the formula \eqref{primepower.convolution} we see that it is equal to the function $u^+$ defined in Definition \ref{uplusdef}. Finally, it can be verified using \eqref{primepower.convolution} that
\[
(p^\alpha)^2 s_1(p^\alpha)*\lambda(p^\alpha) = \sum_{\beta=0}^\alpha p^{2\beta} s_2(p^\beta) \lambda(p^{\alpha-\beta}) = (p^\alpha)^2 s_1^+(p^\alpha),
\]
where $s_1^+$ is defined in Definition \ref{s1plusdef}; therefore the multiplicative function $N^2s_1(N)*\lambda(N)$ is equal to $N^2s_1^+(N)$. This establishes Theorem \ref{g1plus.formula.theorem}.

The proof of Theorem \ref{g1star.formula.theorem} combines the techniques of the above proof with the proof of Theorem \ref{g0star.formula.theorem}, using as a starting point the subspace $\Skstar1N$ of $\Sk1N$ defined by
\begin{equation*}
\Skstar1N = \bigoplus_{d\mid N} i_{1,d,N}\big(\Skplus1n\big),
\end{equation*}
whose dimension can be interpreted as the number of nonisomorphic automorphic representations associated with $\Sk1N$. We omit the details, as by now the method has been amply illustrated.
\end{proof}

\section{Explicit bounds}
\label{explicit.section}

We begin this section by using the formula in Theorem \ref{g0plus.formula.theorem} to extract some explicit bounds on the function $g_0^+(N,2)$, culminating in a proof of Theorem \ref{bennett.bound.theorem}. In the following lemmas, we prove that Theorem \ref{bennett.bound.theorem} holds for certain conveniently chosen classes of integers $N$, after which we combine the results of these lemmas with a modest finite calculation to prove the theorem.

%:Lemma - values of g_0^+ on primes
\begin{lemma}
For every prime $p$, we have $g_0^+(p,2) \le \frac{p+1}{12}$, with equality if and only if $p\equiv11\mod{12}$.
\label{values.on.primes.lemma}
\end{lemma}

\begin{proof}
We directly verify the claim for $p=2$ and $p=3$, so that we may assume $p\ge5$. From Theorem \ref{g0plus.formula.theorem} applied with $k=2$, we have
\[
g_0^+(p) =\tfrac1{12}(p-1) + \begin{cases}
\frac12, &\text{if }p\equiv3\mod4 \\
0, &\text{if }p\equiv1\mod4 \\
\end{cases}\!\!\Bigg\} + \begin{cases}
\frac23, &\text{if }p\equiv2\mod3 \\
0, &\text{if }p\equiv1\mod3 \\
\end{cases}\!\!\Bigg\} - 1.
\]
This establishes the corollary and in fact more, namely that $g_0^+(p)-p/12$ is a constant depending only on the residue class of $p\mod{12}$.
\end{proof}

%:Lemma - bounds for the \nu_i
\begin{lemma}
We have $Ns_0^+(N) \le \phi(N)$, $|\nu^+_2(N)|\le 2^{\omega(N)}$, $|\nu^+_3(N)|\le 2^{\omega(N)}$, and $0\le \nu^+_\infty(N)\le \sqrt N$ for all positive integers $N$.
\label{m.bounds.lemma}
\end{lemma}

\begin{proof}
Since all terms in these four inequalities are multiplicative functions, the asserted inequalities can be checked on prime powers directly from the definitions \ref{s0plusdef}--\ref{nu3plusdef} of the functions $s_0^+$ and $\nu_i^+$. We omit the straightforward verifications.
\end{proof}

%:Corollary - simpler g_0^+ upper bound
\begin{cor}
We have $g_0^+(N,2) \le \tfrac1{12}\phi(N) + \tfrac7{12}2^{\omega(N)} +1$.
\label{g0plus.upper.bound.corollary}
\end{cor}

\begin{proof}
This follows directly from Theorem \ref{g0plus.formula.theorem} and the bounds given in Lemma \ref{m.bounds.lemma}, together with the fact that $|\mu(N)|\le1$.
\end{proof}

%:Lemma - two prime factors
\begin{lemma}
Suppose that $N$ is a composite number with at most two distinct prime factors. Then $g_0^+(N,2) \le \frac{N+1}{12}$, with equality if and only if $N=35$.
\label{two.prime.factors.lemma}
\end{lemma}

\begin{proof}
Since $N$ is composite, it has a divisor $1<d\le\sqrt N$. There are $N/d \ge\sqrt N$ multiples of $d$ less than $N$, none of which is relatively prime to $N$, and hence we have the inequality $\phi(N) \le N-\sqrt N$. From Corollary \ref{g0plus.upper.bound.corollary} and the assumption that $\omega(N)\le2$, we then have
\[
g_0^+(N,2) \le \tfrac1{12}(N-\sqrt N) + \tfrac7{12}2^2 + 1 = \tfrac{N+1}{12} - \big( \tfrac1{12}\sqrt N - \tfrac{13}4 \big).
\]
The quantity $(\tfrac1{12}\sqrt N - \tfrac{13}4)$ is positive as soon as $N>1521$, and so the lemma holds for these large $N$. A direct calculation of $g_0^+(N,2)$ for $N\le1521$ (which discovers the case of equality $N=35$) then shows that the lemma holds for these small $N$ as well.
\end{proof}

%:Lemma - high prime power
\begin{lemma}
Suppose that $N$ is divisible by the sixth power of a prime. Then $g_0^+(N,2) \le \frac{N-6}{12}$.
\label{high.prime.power.lemma}
\end{lemma}

\begin{proof}
Suppose that $p_0^{\alpha_0}$ is a prime power divisor of $N$ with $\alpha_0\ge6$. Then
\[
\frac N{2^{\omega(N)}} = \prod_{p^r\|N} \frac{p^r}2 \ge \frac{p_0^{\alpha_0}}2 \ge \frac{2^{\alpha_0-1}p_0}2 \ge 16p_0,
\]
which is the same as $N/p_0 \ge 16\cdot 2^{\omega(N)}$. Noting that $\phi(N) = N\prod_{p\mid N} (1-\frac1p) \le N (1-\frac1{p_0})$, this implies that
\begin{align*}
N-6 = N\big( 1-\tfrac1{p_0} \big) + \tfrac N{p_0} - 6 &\ge \phi(N) + 16\cdot 2^{\omega(N)} - 6 \\
&\ge \phi(N) + 7\cdot 2^{\omega(N)} + 9\cdot 2^1 - 6 \ge \phi(N) + 7\cdot 2^{\omega(N)} + 12.
\end{align*}
Dividing both sides by 12 and invoking Corollary \ref{g0plus.upper.bound.corollary} establishes the lemma.
\end{proof}

%:Lemma - three prime factors
\begin{lemma}
Suppose that $N$ has at least three distinct prime factors, two of which exceed 5. Then $g_0^+(N,2) \le \frac{N-9}{12}$.
\label{three.prime.factors.lemma}
\end{lemma}

\begin{proof}
Suppose that $p_0<p_1<p_2$ are three distinct prime factors of $N$ with $p_1>5$, so that $p_1\ge7$ and $p_2\ge11$. Then
\[
\frac N{2^{\omega(N)}} \ge \prod_{p\mid N} \frac p2 \ge \frac{p_2}2\frac{p_1}2\frac{p_0}2 \ge \frac{77p_0}8,
\]
which is the same as $\frac N{p_0} \ge \frac{77}8\cdot 2^{\omega(N)}$. This implies that
\begin{align*}
N-9 = N\big( 1-\tfrac1{p_0} \big) + \tfrac N{p_0} - 9 &\ge \phi(N) + \tfrac{77}8\cdot 2^{\omega(N)} - 9 \\
&\ge \phi(N) + 7\cdot 2^{\omega(N)} + \tfrac{21}8 2^3 - 9 \ge \phi(N) + 7\cdot 2^{\omega(N)} + 12
\end{align*}
since $\omega(N)\ge3$. Dividing both sides by 12 and invoking Corollary \ref{g0plus.upper.bound.corollary} establishes the lemma.
\end{proof}

%:Lemma - 2 or 3 divides N plus a big prime
\begin{lemma}
If $(N,6)>1$ and $N$ has a prime factor exceeding 41, then $g_0^+(N,2) \le \frac N{12}$.
\label{2or3.divides.n.lemma}
\end{lemma}

\begin{proof}
Since either $2\mid N$ or $3\mid N$, we have $\phi(N)\le\frac{2N}3$. Let $p>41$ be a prime factor of $N$. Then
\[
\frac N{2^{\omega(N)}} \ge \prod_{p\mid N} \frac p2 \ge \frac{43}2,
\]
which is the same as $\frac7{12}\cdot 2^{\omega(N)} \le \frac{7N}{258}$. Then by Corollary \ref{g0plus.upper.bound.corollary},
\[
g_0^+(N,2) \le \tfrac1{12}\phi(N) + \tfrac7{12}2^{\omega(N)} + 1 \le \tfrac1{12}\tfrac{2N}3 + \tfrac{7N}{258} + 1 = \tfrac N{12} - \big( \tfrac N{1548} - 1 \big).
\]
This establishes the lemma for $N\ge1548$, and we check by direct calculation that the lemma holds for $N<1548$.
\end{proof}

%:Proof of Bennett's bound
\begin{proof}[Proof of Theorem \ref{bennett.bound.theorem}]
Lemmas \ref{values.on.primes.lemma} and \ref{two.prime.factors.lemma} show that if $\omega(N)\le 2$, then $g_0^+(N,2) \le \frac{N+1}{12}$ with equality if and only if either $N=35$ or $N$ is a prime that is congruent to $11\mod{12}$. It remains to show that $g_0^+(N,2) < \frac{N+1}{12}$ when $\omega(N)\ge3$. This inequality follows from Lemma \ref{three.prime.factors.lemma} if two of the prime factors of $N$ exceed 5; therefore we need only consider numbers of the form $N=2^{\alpha_1}3^{\alpha_2}5^{\alpha_3}p^{\alpha_4}$ with $p>5$ and at least three of the $\alpha_i$ positive. No such integer can be relatively prime to 6, however; thus if $p>41$, we have $g_0^+(N,2) < \frac{N+1}{12}$ by Lemma \ref{2or3.divides.n.lemma}. Furthermore, if any $\alpha_i\ge6$, then $g_0^+(N,2) < \frac{N+1}{12}$ by Lemma \ref{high.prime.power.lemma}.

Consequently, the only integers $N$ for which we have not verified the theorem are those of the form $N=2^{\alpha_1}3^{\alpha_2}5^{\alpha_3}p^{\alpha_4}$ with $7\le p\le 41$, where each $0\le\alpha_i\le 5$ and at least three of the $\alpha_i$ are positive. There are 10,125 integers of this form, and a direct calculation verifies that $g_0^+(N,2) \le \frac N{12}-\frac32$ for these integers. This establishes the theorem.
\end{proof}

We now turn to the evaluation of $g_0^+(2^\alpha N,k)$.

%:Proof of Bennett's inequality
\begin{proof}[Proof of Theorem \ref{bennett.inequality.theorem}]
Let $N\ge3$ be an odd squarefree integer, and let $\alpha\ge4$ be an integer. Then from Theorem \ref{g0plus.formula.theorem},
\begin{multline}
g_0^+(2^\alpha N) = \tfrac{k-1}{12}2^\alpha N s_0^+(2^\alpha)s_0^+(N) - \tfrac12 \nu_\infty^+(2^\alpha)\nu_\infty^+(N) \\ {}+ c_2(k) \nu_2^+(2^\alpha)\nu_2^+(N) + c_3(k) \nu_3^+(2^\alpha)\nu_3^+(N) + \delta(\tfrac k2) \mu(2^\alpha) \mu(N).
\label{g0plus.2alphaN.formula}
\end{multline}
From Definitions \ref{nuinftyplusdef}--\ref{nu3plusdef}, we see that $\nu_\infty^+(2^\alpha) = \nu_2^+(2^\alpha) = \nu_3^+(2^\alpha) = \mu(2^\alpha) = 0$ since $\alpha\ge4$ and $N\ge3$ is squarefree. Also, from Definition \ref{s0plusdef},
\[
s_0^+(2^\alpha) s_0^+(N) = \big( 1-\tfrac12 \big) \big( 1-\tfrac1{2^2} \big) \prod_{p\mid N} \big( 1-\tfrac1p \big) = \tfrac38\tfrac{\phi(N)}N.
\]
We conclude that $g_0^+(2^\alpha N) = \tfrac{k-1}{12}2^\alpha N \cdot \tfrac38\tfrac{\phi(N)}N = (k-1)2^{\alpha-4}\phi(N)$ as asserted.

Now considering equation \eqref{g0plus.2alphaN.formula} in the case $\alpha=1$, we have
\begin{align*}
g_0^+(2N) &= \tfrac{k-1}{12}2N s_0^+(2)s_0^+(N) - \tfrac12 \nu_\infty^+(2)\nu_\infty^+(N) \\
&\qquad{}+ c_2(k) \nu_2^+(2)\nu_2^+(N) + c_3(k) \nu_3^+(2)\nu_3^+(N) + \delta(\tfrac k2) \mu(2) \mu(N) \\
&= \tfrac{k-1}{12} \phi(N) - c_2(k)\nu_2^+(N) -2c_3(k) \nu_3^+(N) - \delta(\tfrac k2)\mu(N).
\end{align*}
By the bounds given in Lemma \ref{m.bounds.lemma} and the definitions \ref{c2def}--\ref{c3def} of $c_2(k)$ and $c_3(k)$,
\[
g_0^+(2N) \le \tfrac{k-1}{12}\phi(N) + \tfrac14 2^{\omega(N)} + \tfrac23 2^{\omega(N)} + 1 = \tfrac{k-1}{12}\phi(N) + \tfrac{11}{12}2^{\omega(N)} + 1.
\]
Since $\phi(N) = \prod_{p\mid N} (p-1)$ and $2^{\omega(N)} = \prod_{p\mid N} 2$, we have $2^{\omega(N)}\le \phi(N)$ with equality if and only if $N=3$. We verify by hand that $g_0^+(6,k) \le 2(k-1)$, which takes care of the case $N=3$. When $N>3$, we have
\[
g_0^+(2N) < \tfrac{k-1}{12}\phi(N) + \tfrac{11}{12}\phi(N) + 1 \le (k-1)\phi(N) + 1,
\]
which establishes the last claim of the theorem.
\end{proof}

\section{Calculations of values of $g_0^+(N,2)$ and $g_0(N,2)$}
\label{up.to.100.section}

Using the formula given in Theorem \ref{g0plus.formula.theorem}, we can derive explicit inequalities for the function $g_0^+(N,2)$. We can thus determine the precise preimage of any fixed value of $g_0^+(N,2)$ by combining these inequalities with finite computations. We remark that Halberstadt and Kraus \cite{HK} independently employed similar methods in their calculations of the set of integers for which $g_0^+(N,2)\le3$.

We begin by stating a few lemmas giving simple but explicit inequalities for the multiplicative functions that concern us. We remind the reader of the definition \eqref{A0plus.def} of the constant $A_0^+$:
\[\textstyle
A_0^+ = \prod_p \big( 1-\frac1{p^2-p} \big) \approx 0.373956.
\]

%:Lemma - lower bound for Ns_0^+(N)
\begin{lemma}
We have $Ns_0^+(N) > A_0^+\phi(N)$ for all integers $N\ge1$.
\label{Ns0plus.lower.bound.lemma}
\end{lemma}

\begin{proof}
From the definition \ref{s0plusdef} of $s_0^+$, we see that on prime powers
\[
p^rs_0^+(p^r) \ge p^r\big( 1-\tfrac1p-\tfrac1{p^2} \big) = \tfrac{p^{r-1}(p^2-p-1)(p-1)}{p(p-1)} = \phi(p^r) \big( 1-\tfrac1{p^2-p} \big).
\]
Therefore
\[
Ns_0^+(N) = \prod_{p^r\|N} p^rs_0^+(p^r) \ge \prod_{p^r\|N} \phi(p^r) \prod_{p\mid N} \big( 1-\tfrac1{p^2-p} \big) > \phi(N)\cdot A_0^+
\]
as claimed.
\end{proof}

%:Lemma - omega upper bound
\begin{lemma}
We have $2^{\omega(N)} \le 2^{4-\frac{\log16}{\log11}} N^{\frac{\log2}{\log11}}$ for all $N\ge1$.
\label{omega.upper.bound.lemma}
\end{lemma}

\begin{proof}
We have
\begin{align*}
2^{\omega(N)} = \bigg( \prod_{\substack{p\mid N \\ p\le7}} 2 \bigg) \bigg( \prod_{\substack{p\mid N \\ p\ge11}} 2 \bigg) &\le \bigg( \prod_{\substack{p\mid N \\ p\le7}} 2 \big( \tfrac p2 \big)^{\frac{\log 2}{\log 11}} \bigg) \bigg( \prod_{\substack{p\mid N \\ p\ge11}} p^{\frac{\log 2}{\log 11}} \bigg) \\
&\le \bigg( \prod_{\substack{p\mid N \\ p\le7}} 2^{1-\frac{\log 2}{\log 11}} \bigg) \bigg( \prod_{p\mid N} p^{\frac{\log 2}{\log 11}} \bigg) \le 2^{4(1-\frac{\log 2}{\log 11})} N^{\frac{\log 2}{\log 11}}
\end{align*}
as claimed.
\end{proof}

%:Lemma - phi lower bound
\begin{lemma}
We have $\phi(N) \ge \frac{N\log 2}{\log2N}$ for $N\ge2$.
\label{phi.lower.bound.lemma}
\end{lemma}

\begin{proof}
This is Theorem 3.1(g) of Bressoud and Wagon \cite{BW}.
\end{proof}

%:Proposition - g_0^+ greater than 100
\begin{prop}
We have $g_0^+(N,2)>100$ for all $N>132\text{\rm,}000$.
\label{greater.than.100.prop}
\end{prop}

\begin{proof}
Suppose first that $N$ is not a perfect square. Then $\nu^+_\infty(N)=0$ by Definition \ref{nuinftyplusdef}, while $c_2(2)=-\frac14$ and $c_3(2)=-\frac13$ by Definitions \ref{c2def}--\ref{c3def}. Therefore the formula in Theorem \ref{g0plus.formula.theorem}, applied with $k=2$, implies the inequality
\[
g_0^+(N,2) \ge \tfrac1{12}Ns_0^+(N) -\tfrac14 |\nu^+_2(N)| -\tfrac13 |\nu^+_3(N)| - \big| \delta\big( \tfrac k2 \big)\mu(N) \big|.
\]
Applying Lemmas \ref{m.bounds.lemma} and \ref{Ns0plus.lower.bound.lemma}, and noting that $|\delta(\tfrac k2)\mu(N)|\le1$, gives
\begin{equation}
g_0^+(N,2) > \tfrac{A_0^+}{12}\phi(N) - \tfrac7{12}2^{\omega(N)} - 1.
\label{g0plus.lower.bound}
\end{equation}
From Lemmas \ref{omega.upper.bound.lemma} and \ref{phi.lower.bound.lemma} we conclude that
\[
g_0^+(N,2) > \tfrac{A_0^+N\log2}{12\log2N} - \tfrac7{12}2^{4-\frac{\log16}{\log11}} N^{\frac{\log2}{\log11}} - 1.
\]
It can be verified that the right-hand side is an increasing function of $N$ for $N>9$,000 and takes a value exceeding 100 when $N=1$32,000. This establishes the theorem in the case where $N$ is not a perfect square.

Suppose now that $N=M^2$ is a perfect square, where $M\ge1$. Then the formula in Theorem \ref{g0plus.formula.theorem}, applied with $k=2$, implies
\[
g_0^+(M^2,2) \ge \tfrac1{12}M^2s_0^+(M^2) -\tfrac12\nu_\infty^+(M^2) -\tfrac14 |\nu^+_2(M^2)| -\tfrac13 |\nu^+_3(M^2)|
\]
since $\mu(M^2)=0$. Applying Lemmas \ref{m.bounds.lemma} and \ref{Ns0plus.lower.bound.lemma} gives
\begin{equation}
g_0^+(M^2,2) > \tfrac{A_0^+}{12}\phi(M^2) - \tfrac12\sqrt{M^2} - \tfrac7{12}2^{\omega(M^2)} = \tfrac{A_0^+}{12}M\phi(M) - \tfrac12M - \tfrac7{12}2^{\omega(M)},
\label{g0plus.square.lower.bound}
\end{equation}
using the elementary facts that $\phi(M^2)=M\phi(M)$ and $\omega(M^2)=\omega(M)$. From Lemmas \ref{omega.upper.bound.lemma} and \ref{phi.lower.bound.lemma} we conclude that
\[
g_0^+(m^2,2) > \tfrac{A_0^+M^2\log2}{12\log2M} - \tfrac12M - \tfrac7{12}2^{4-\frac{\log16}{\log11}} M^{\frac{\log2}{\log11}} - 1.
\]
It can be verified that the right-hand side is an increasing function of $M$ for $M>170$ and takes a value exceeding 100 when $M=280$. Since $280^2=\text{78,400}<\text{132,000}$, this establishes the theorem in this case as well.
\end{proof}

Using the formula in Theorem \ref{g0plus.formula.theorem}, it takes only a couple of minutes to compute $g_0^+(N,2)$ for all $N\le1$32,000. We discover that there are exactly 2,965 integers $N$ for which $g_0^+(N,2)\le100$. For example, there are exactly 40 solutions to the equation $g_0^+(N,2)=100$, namely
\begin{multline*}
N = 1213, 1331, 2169, 2583, 2662, 2745, 3208, 3232, 3465, 3608, 4040, 4302, 4338, \\ 4772, 4804, 4848, 5084, 5092, 5166, 5252, 5324, 5490, 5572, 5904, 6336, 6820, 6930, \\ 7056, 7188, 7212, 7920, 8052, 8484, 8652, 8676, 8940, 9060, 10332, 10980, 13860.
\end{multline*}
We found that for every integer $0\le k\le100$ there are at least 13 solutions to the equation $g_0^+(N,2)=k$, and there are only 13 solutions for $k=86$. The largest number of solutions for $k$ in this range is 68, attained by $k=96$.

As $N$ ranges from 1 to 132,000, the values taken by $g_0^+(N,2)$ include every nonnegative integer up to and including 4,361. In total, we found 9,566 of the integers less than 10,000 among the values of $g_0^+(N,2)$ during this calculation, and of course extending the range of computation further would likely increase this number. The following assertion therefore seems reasonable:

%:Conjecture - g_0^+(N,2) is surjective
\begin{conjecture}
For every nonnegative integer $k$, there is at least one positive integer $N$ such that $g_0^+(N,2)=k$.
\label{g0plus.surjective.conjecture}
\end{conjecture}

However, we can show that the analogous conjecture is false for $g_0(N,2)$. The results of our computations are as follows:

%:Proposition - g_0(N,2) not surjective
\begin{prop}
The equation $g_0(N,2)=k$ has a solution $N$ for every integer $0\le k\le 1000$ except  for {\rm $k=150$, 180, 210, 286, 304, 312, 336, 338, 348, 350, 480, 536, 570, 598, 606, 620, 666, 678, 706, 730, 756, 780, 798, 850, 876, 896, 906, 916, and 970.}
\label{not.surjective.prop}
\end{prop}

\begin{proof}
In analogy with Lemma \ref{m.bounds.lemma}, the inequalities
\begin{equation}
0\le\nu_2(N)\le 2^{\omega(N)}, \quad 0\le\nu_3(N)\le 2^{\omega(N)}, \quad\text{and } 0\le \nu_\infty(N)\le \sqrt N s_0(N)
\label{analogy.inequalities}
\end{equation}
follow easily by considering the values of all expressions involved on prime powers. Using these inequalities, Proposition \ref{g0.formula.prop} provides the lower bound
\begin{align*}
g_0(N,2) &= \tfrac1{12}Ns_0(N) - \tfrac12\nu_\infty(N) - \tfrac14\nu_2(N) - \tfrac13\nu_3(N) + 1 \\
&> \tfrac1{12}Ns_0(N) - \tfrac12\sqrt Ns_0(N) - \tfrac14 2^{\omega(N)} - \tfrac13 2^{\omega(N)} \\
&> \tfrac1{12}(N-6\sqrt N)s_0(N) - \tfrac7{12} 2^{\omega(N)}.
\end{align*}
If $N>36$ then $N-6\sqrt N>0$, and so we can use the bound $s_0(N)\ge1$ (which follows directly from Definition \ref{s0def}) and Lemma \ref{omega.upper.bound.lemma} to obtain the inequality
\[
g_0(N,2) > \tfrac1{12}(N-6\sqrt N) - \tfrac7{12} 2^{4-\frac{\log16}{\log11}} N^{\frac{\log2}{\log11}}.
\]
It is easily shown from this inequality that if $N>1$3,500, then $g_0(N,2) > 1$,000. A calculation of all of the values of $g_0(N,2)$ as $N$ ranges up to 13,500 shows that the 29 integers listed in the statement of the proposition are not in fact values of $g_0(N,2)$, while the other 972 integers between 0 and 1,000 are.
\end{proof}

We have extended these computations for $N$ ranging up to 124,000; it can be shown in a manner similar to the proof of Proposition \ref{not.surjective.prop} that this is guaranteed to find all solutions to $g_0(N,2) \le 1$0,000. Based on this numerical evidence, it seems that approximately 94-95\% of all positive integers are values of $g_0(N,2)$. However, because $g_0(N,2)$ is not a multiplicative function but rather a linear combination of multiplicative functions, we do not know how to approach the problem of determining the density of its range. In particular, we cannot prove that a positive proportion of integers are omitted as values (as the data leads us to suspect); indeed, we cannot even prove that a positive proportion of integers are {\em taken} as values of $g_0(N,2)$.

Certainly, there do not seem to be any residue classes of values that are systematically omitted by the function $g_0(N,2)$, so a proof that a positive proportion of integers are omitted seems nontrivial. In fact, these values seem to be quite well distributed among residue classes in general. One notable exception is that $g_0(N,2)-1$ is noticeably more likely to be divisible by powers of 2 then random integers; this is not surprising in hindsight, since the multiplicative functions involved in the formula in Proposition \ref{g0.formula.prop} all have the tendency to take even values on prime powers. Every odd integer below 10,000 is taken as a value of $g_0(N,2)$, but we do not know whether or not this trend persists.

\section{Minimal and maximal orders}
\label{min.max.orders.section}

In this section we provide the arguments necessary to convert the exact formulas for $g_0$, $g_0^*$, $g_0^+$, $\rho_0$, and $g_1$ and its variants into asymptotic upper and lower bounds. We remark again that each of these bounds is sharp, and the avid reader wil be able to convert the proofs below into constructions of sequences of integers that attain the bounds in question. We begin with three simple lemmas concerning the order of growth of some of the multiplicative functions we have encountered.

%:Lemma - converging to 6/\pi^2
\begin{lemma}
We have
\[
\prod_{p\le y} \big( 1-\tfrac1{p^2} \big) = \tfrac6{\pi^2} + O\big( \tfrac1y \big)
\]
for all $y\ge2$.
\label{converge.to.lemma}
\end{lemma}

\begin{proof}
The product in question converges to $\prod_p (1-\frac1{p^2}) = \frac1{\zeta(2)}$ as $y$ tends to infinity. To assess the error term for the partial product, note that
\[
\sum_{p>y} \log\big( 1-\tfrac1{p^2} \big)^{-1} \ll \sum_{p>y} \tfrac1{p^2} < \sum_{n>y} \tfrac1{n^2} \ll \tfrac1y.
\]
Therefore
\[
\prod_{p>y} \big( 1-\tfrac1{p^2} \big)^{-1} = \exp\!\big( O\big( \tfrac1y \big) \big) = 1 + O\big( \tfrac1y \big),
\]
which implies that
\[
\prod_{p\le y} \big( 1-\tfrac1{p^2} \big) = \prod_p \big( 1-\tfrac1{p^2} \big) \prod_{p>y} \big( 1-\tfrac1{p^2} \big)^{-1} = \tfrac1{\zeta(2)} \big( 1 + O\big( \tfrac1y \big) \big) = \tfrac6{\pi^2} + O\big( \tfrac1y \big),
\]
since $\zeta(2)=\tfrac{\pi^2}6$.
\end{proof}

%:Lemma - s0 and log log
\begin{lemma}
We have
\[
1 \le s_0(N) \le \tfrac{6e^\gamma}{\pi^2}\log\log N + O(1)
\]
uniformly for all integers $N\ge2$.
\label{s0.and.log.log.lemma}
\end{lemma}

\begin{proof}
The lower bound $s_0(N)\ge1$ is trivial. For the upper bound, first we consider the special case where $N$ has the form $N_y = \prod_{p\le y} p$. In this case,
\[
s_0(N_y) = \prod_{p\le y} \big( 1+\tfrac1p \big) = \prod_{p\le y} \big( 1-\tfrac1p \big)^{-1} \prod_{p\le y} \big( 1-\tfrac1{p^2} \big).
\]
An asymptotic formula for the first product on the right-hand side is well known: Mertens' formula is
\[
\prod_{p\le y} \big( 1-\tfrac1p \big)^{-1} = e^\gamma\log y + O(1).
\]
Therefore
\[
s_0(N_y) = \big( e^\gamma\log y + O(1) \big) \big( \tfrac6{\pi^2} + O\big( \tfrac1y \big) \big) = \tfrac{6e^\gamma}{\pi^2}\log y + O(1)
\]
by Lemma \ref{converge.to.lemma}. On the other hand, the prime number theorem tells us that
\[
\log N_y = \sum_{p\le y} \log p = y \big( 1 + O\big( \tfrac 1{\log y} \big) \big)
\]
(in fact we could be much more generous with the error term if need be). Therefore
\[
s_0(N_y) = \tfrac{6e^\gamma}{\pi^2}\log\log N_y + O(1),
\]
which establishes the lemma for integers of the form $N_y$.

Now consider an arbitrary integer $N\ge2$. Choose $y$ to be the $\omega(N)$th prime number, and set $N_y = \prod_{p\le y} p$ as before. Then $N\ge N_y$, and the various prime factors of $N$ are at least as large as the corresponding prime factors of $N_y$. Therefore
\[
s_0(N) = \prod_{p\mid N} \big( 1+\tfrac1p \big) \le \prod_{p\le y} \big( 1+\tfrac1p \big) = \tfrac{6e^\gamma}{\pi^2}\log\log N_y + O(1) \le \tfrac{6e^\gamma}{\pi^2}\log\log N + O(1)
\]
as desired.
\end{proof}

%:Lemma - u(N) and t(N)
\begin{lemma}
We have $t(N) \le u(N) \le N\tau(N)$ for all $N\ge1$.
\label{un.and.tn.lemma}
\end{lemma}

\begin{proof}
Since all three functions are multiplicative and nonnegative, it suffices to show that $t(p^\alpha) \le u(p^\alpha) \le p^\alpha\tau(p^\alpha)$ for all prime powers $p^\alpha$. This is easily verified by hand when $\alpha=1$ and $\alpha=2$. When $\alpha\ge3$, we need to show that
\[
p^{\alpha-4}(p-1)^3((\alpha+1)p+\alpha-3) \le p^{\alpha-2}(p-1)((\alpha+1)p+\alpha-1) \le p^{\alpha}(\alpha+1)
\]
for all primes $p\ge2$. The first inequality follows from the obvious inequality \[(p-1)^2((\alpha+1)p+\alpha-3) \le p^2((\alpha+1)p+\alpha-1)\] upon multiplying through by $p^{\alpha-4}(p-1)$, and the second inequality similarly follows from \[(p-1)((\alpha+1)p+\alpha-1) \le p((\alpha+1)p+\alpha+1)\] upon multiplying through by $p^{\alpha-2}$.
\end{proof}

\begin{proof}[Proof of Theorem \ref{asymptotic.bounds.g0.theorem}]
Starting with the formula
\[
g_0(N,k) = \tfrac{k-1}{12}Ns_0(N) - \tfrac12\nu_\infty(N)+c_2(k) \nu_2(N) +c_3(k) \nu_3(N) + \delta\big( \tfrac k2 \big)
\]
given by Proposition \ref{g0.formula.prop}, we use the inequalities \eqref{analogy.inequalities} to deduce that
\begin{multline*}
\tfrac{k-1}{12}Ns_0(N) - \tfrac12\sqrt N\,s_0(N) - |c_2(k)| 2^{\omega(N)} - |c_3(k)| 2^{\omega(N)} \le g_0(N,k) \\ \le \tfrac{k-1}{12}Ns_0(N) + |c_2(k)| 2^{\omega(N)} + |c_3(k)| 2^{\omega(N)} + 1.
\end{multline*}
The coefficients $c_2(k)$ and $c_3(k)$ are uniformly bounded, and $2^{\omega(N)} \ll \sqrt N$. Therefore we may write these inequalities as
\[
\tfrac{k-1}{12}Ns_0(N) + O(\sqrt N\,s_0(N)) \le g_0(N,k) \le \tfrac{k-1}{12}Ns_0(N) + O(2^{\omega(N)}).
\]
By Lemma \ref{s0.and.log.log.lemma}, we conclude that
\[
\tfrac{k-1}{12}N + O(\sqrt N\log\log N) \le g_0(N,k) \le \tfrac{k-1}{12}N \big( \tfrac{6e^\gamma}{\pi^2}\log\log N + O(1) \big) + O(2^{\omega(N)}),
\]
which establishes Theorem \ref{asymptotic.bounds.g0.theorem}(a).

In a similar way, combining the formula
\[
g_0^*(N,k) = \tfrac{k-1}{12}Ns_0^*(N) - \tfrac12\nu^*_\infty(N)+c_2(k) \nu^*_2(N) +c_3(k) \nu^*_3(N) + \delta\big( \tfrac k2 \big)\delta(N),
\]
from Theorem \ref{g0star.formula.theorem} with the easily verifiable inequalities
\[
\tfrac6{\pi^2} = \tfrac1{\zeta(2)} < s_0^*(N) \le 1,\, |\nu_2^*(N)|\le1,\, |\nu_3^*(N)|\le1,\, \text{ and } 0\le \nu_\infty^*(N)\le \tfrac{\phi(N)}{\sqrt N}
\]
establishes Theorem \ref{asymptotic.bounds.g0.theorem}(b). Moreover, combining the formula
\[
g_0^+(N,k) = \tfrac{k-1}{12}Ns_0^+(N) - \tfrac12\nu^+_\infty(N)+c_2(k) \nu^+_2(N) +c_3(k) \nu^+_3(N) + \delta\big( \tfrac k2 \big)\delta(N),
\]
from Theorem \ref{g0star.formula.theorem} with the inequalities from Lemma \ref{m.bounds.lemma} and the additional inequality
\[
Ns_0^+(N) \ge \phi(N) \prod_{p\mid N} \big( 1-\tfrac1p \big)^{-1} \big( 1-\tfrac1p -\tfrac1{p^2} \big) = \prod_{p\mid N} \big( 1-\tfrac1{p^2-p} \big) > A_0^+,
\]
which follows from the definition \eqref{A0plus.def} of $A_0^+$, establishes Theorem \ref{asymptotic.bounds.g0.theorem}(c).
\end{proof}

The proof of Theorem \ref{asymptotic.bounds.g1.theorem} is very similar, and we omit the details except to mention that Lemma \ref{un.and.tn.lemma} plays a role in simplifying the error terms. As for Theorem \ref{asymptotic.bounds.rhos.theorem}, we can investigate the size of $\rho_0(N,k)$ (for example) using the information discovered in the proof of Theorem \ref{asymptotic.bounds.g0.theorem}. We saw that
\[
g_0^+(N,k) = \tfrac{k-1}{12}Ns_0^+(N) + O(\sqrt N) = \tfrac{k-1}{12}Ns_0^+(N) \big( 1 + O\big( \tfrac{\log\log N}{\sqrt N} \big) \big)
\]
and similarly $g_0(N,k) = \tfrac{k-1}{12}Ns_0(N) (1 + O(\tfrac{\log\log N}{\sqrt N}))$. Therefore when $g_0(N,k)\ne0$, we have
\[
\rho_0(N,k) = \tfrac{g_0^+(N,k)}{g_0(N,k)} = \tfrac{s_0^+(N)}{s_0(N)} \big( 1 + O\big( \tfrac{\log\log N}{\sqrt N} \big) \big).
\]
The size of the multiplicative function $\tfrac{s_0^+(N)}{s_0(N)}$ can be investigated as in the proof of Lemma \ref{s0.and.log.log.lemma}. We find that
\[
\tfrac{A_0^+\pi^2}{6e^{2\gamma}(\log\log N)^2} \big( 1 + O\big( \tfrac1{\log\log N} \big) \big) < \tfrac{s_0^+(N)}{s_0(N)} \le 1,
\]
which is enough to establish Theorem \ref{asymptotic.bounds.rhos.theorem}(a). The proof of Theorem \ref{asymptotic.bounds.rhos.theorem}(b) is quite similar.

\section{Average orders}
\label{average.section}

In this final section we prove Theorems \ref{average.order.g0.theorem}, \ref{average.order.g1.theorem}, and \ref{average.order.rhos.theorem}. As it happens, the multiplicative functions under consideration are all in a class of multiplicative functions whose average orders can be calculated rather easily. The following proposition is representative of the average-order theorems for multiplicative functions in the literature; we include a proof for the sake of completeness.

%:Proposition - general average orders
\begin{prop}
Suppose that $h(n)$ is a multiplicative function with the property that for some positive constant $\eta$, we have $(h*\mu)(n) \ll n^{-\eta}$. Then for any $\beta>-1$, we have
\[
\sum_{n\le x} n^\beta h(n) \sim \frac{c(h)x^{\beta+1}}{\beta+1},
\]
where
\[
c(h) = \prod_p \Big( 1-\frac1p \Big) \Big( 1 + \frac{h(p)}p + \frac{h(p^2)}{p^2} + \cdots \Big).
\]
In particular, the average order of the function $n^\beta h(n)$ is $c(h)n^\beta$.
\label{general.average.orders.prop}
\end{prop}

\begin{proof}
Let $g$ denote the convolution $h*\mu$, so that $h(n) = \sum_{d\mid n} g(d)$ by the M\"obius inversion formula; we note that $g$ is multiplicative as well. For $x\ge1$ we have
\begin{multline*}
\sum_{n\le x} n^\beta h(n) = \sum_{n\le x} n^\beta \sum_{d\mid n} g(d) = \sum_{d\le x} g(d) \sum_{\substack{n\le x \\ d\mid n}} n^\beta \\
= \sum_{d\le x} g(d) \sum_{md\le x} (dm)^\beta = \sum_{d\le x} d^\beta g(d) \sum_{m\le x/d} m^\beta.
\end{multline*}
Using the fact that
\[
\sum_{m\le y} m^\beta = \frac{y^{\beta+1}}{\beta+1} + O(y^\beta)
\]
for any fixed $\beta>-1$, we see that
\begin{align}
\sum_{n\le x} n^\beta h(n) &= \sum_{d\le x} d^\beta g(d) \Big( \frac{(x/d)^{\beta+1}}{\beta+1} + O\big( (x/d)^\beta \big) \Big) \notag \\
&= \frac{x^{\beta+1}}{\beta+1} \sum_{d\le x} \frac{g(d)}d + O\bigg( x^\beta \sum_{d\le x} |g(d)| \bigg).
\label{into.bits}
\end{align}
Since $g(d) \ll d^{-\eta}$, the sum in the main term is a truncation of a convergent sum, as the tail can be estimated by
\[
\sum_{d>x} \frac{g(d)}d \ll \sum_{d>x} d^{-\eta-1} \ll x^{-\eta}.
\]
Moreover, since $g$ is multiplicative we can write
\begin{equation}
\sum_{n=1}^\infty \frac{g(d)}d = \prod_p \Big( 1 + \frac{g(p)}p + \frac{g(p^2)}{p^2} + \cdots \big).
\label{g.first}
\end{equation}
Since $h(p^\alpha) - h(p^{\alpha-1}) = g(p^\alpha)$, it is easily seen that
\[
\Big( 1-\frac1p \Big) \Big( 1 + \frac{h(p)}p + \frac{h(p^2)}{p^2} + \cdots \Big) = 1 + \frac{g(p)}p + \frac{g(p^2)}{p^2} + \cdots,
\]
where convergence is ensured by the hypothesis $g(p^\alpha) \ll p^{-\eta\alpha}$. Therefore equation \eqref{g.first} becomes
\[
\sum_{d=1}^\infty \frac{g(d)}d = \prod_p \Big( 1-\frac1p \Big) \Big( 1 + \frac{h(p)}p + \frac{h(p^2)}{p^2} + \cdots \Big) = c(h).
\]
Finally, we have the estimate
\[
\sum_{d\le x} |g(d)| \ll \sum_{d\le x} d^{-\eta} \ll E_\eta(x),
\]
where
\[
E_\eta(x) =\begin{cases}
x^{1-\eta}, &\text{if } 0<\eta<1, \\
\log x, &\text{if } \eta=1, \\
1, &\text{if } \eta>1. \\
\end{cases}
\]
Assembling this information and applying it to equation \eqref{into.bits} yields
\begin{align*}
\sum_{n\le x} n^\beta h(n) &= \frac{x^{\beta+1}}{\beta+1} \bigg( \sum_{d=1}^\infty \frac{g(d)}d + O\bigg( \sum_{d>x} \frac{g(d)}d \bigg) \bigg) + O\bigg( x^\beta \sum_{d\le x} |g(d)| \bigg) \\
&= \frac{x^{\beta+1}}{\beta+1} c(h) + O( x^{\beta+1-\eta} + x^\beta E_\eta(x) ) \\
&= \frac{c(h)x^{\beta+1}}{\beta+1} + O( x^\beta E_\eta(x) ),
\label{into.bits}
\end{align*}
which establishes the proposition.
\end{proof}

To apply this proposition to prove Theorem \ref{average.order.g0.theorem}(a), for example, we start with the equation $g_0(N,k) = \tfrac{k-1}{12}Ns_0(N) + O(\sqrt N\log\log N)$. It follows that
\begin{equation}
\sum_{N\le x} g_0(N,k) = \tfrac{k-1}{12} \sum_{N\le x} Ns_0(N) + O(x^{3/2}\log\log x).
\label{only.main.term.equation}
\end{equation}
We note that the function $s_0\ast\mu$ is multiplicative and satisfies $s_0(p)=\frac1p$ and $s_0(p^\alpha)=0$ for $\alpha\ge2$. Therefore the hypothesis of Proposition \ref{general.average.orders.prop} is satisfied with $\eta=1$, and so we conclude that
\[
\sum_{N\le x} Ns_0(N) \sim \tfrac12c(s_0)x^2,
\]
where
\begin{align*}
c(s_0) &= \prod_p \Big( 1-\frac1p \Big) \Big( 1 + \frac{s_0(p)}p + \frac{s_0(p^2)}{p^2} + \cdots \Big) \\
&= \prod_p \Big( 1-\frac1p \Big) \Big( 1 + \Big( 1+\frac1p \Big) \Big( \frac1p + \frac1{p^2} + \cdots \Big) \Big) \\
&= \prod_p \Big( 1+\frac1{p^2} \Big) \\
&= \prod_p \Big( 1-\frac1{p^2} \Big)^{\!-1} \prod_p \Big( 1-\frac1{p^4} \Big)\\
&= \frac{\zeta(2)}{\zeta(4)} = \frac{\pi^2/6}{\pi^4/90} = \frac{15}{\pi^2}.
\end{align*}
Combining this with equation \eqref{only.main.term.equation}, we conclude that
\[
\sum_{N\le x} g_0(N,k) \sim \frac{k-1}{12} \frac{15}{\pi^2} \frac{x^2}2 = \frac{5(k-1)x^2}{8\pi^2},
\]
which implies that the average order of $g_0(N,k)$ is indeed $\frac{5(k-1)N}{4\pi^2}$. The proofs of the other seven average-order assertions in Theorems \ref{average.order.g0.theorem}, \ref{average.order.g1.theorem}, and \ref{average.order.rhos.theorem} all follow this outline, and we omit the details of the calculations.

\bigskip
{\noindent\small{\it Acknowledgements.} The author is grateful to Mike Bennett and Nike Vatsal for sharing their expertise and bringing these problems to his attention. The author acknowledges the support of the Natural Sciences and Engineering Research Council.}

\bibliographystyle{amsplain}
\bibliography{newforms}

\end{document}